\documentclass[colorlinks]{elsarticle}
\usepackage{amsmath}
\usepackage{amsthm}
\usepackage{amssymb}
\usepackage[numbers]{natbib}
\usepackage{tabularx}
\usepackage{float}
\usepackage[ruled,linesnumbered,vlined]{algorithm2e} 
\usepackage[table]{xcolor} 
\usepackage{caption}
\usepackage{subcaption,graphicx}
\usepackage[export]{adjustbox}
\usepackage{lineno,hyperref}
\usepackage{varioref} 
\usepackage{refstyle}
\usepackage{multirow} 
\usepackage [english]{babel}
\usepackage [autostyle, english = american]{csquotes}
\usepackage[super]{nth} 
\MakeOuterQuote{"}
\allowdisplaybreaks

\modulolinenumbers[5]

\usepackage[colorinlistoftodos]{todonotes}

\usepackage{changes}

\usepackage{siunitx}
\usepackage{tikz}
\usepackage{pgfplots}
\pgfplotsset{width=10cm,compat=1.9}
\definecolor{pistachio}{rgb}{0.58, 0.77, 0.45}
\definecolor{cinnamon}{rgb}{0.82, 0.41, 0.12}
\definecolor{celadon}{rgb}{0.67, 0.88, 0.69}
\definecolor{cadmiumgreen}{rgb}{0.0, 0.42, 0.24}
\definecolor{amethyst}{rgb}{0.6, 0.4, 0.8}
\definechangesauthor[name=Ruslan Krenzler,color=cadmiumgreen]{RK}
\definechangesauthor[name=Lin Xie,color=amethyst]{LX}
\definechangesauthor[name=Nils Thieme, color=amethyst]{NT}

\theoremstyle{plain}

\theoremstyle{remark}

\theoremstyle{definition}

\journal{EJOR}

\biboptions{authoryear}

\begin{document}

\begin{frontmatter}

\title{Formulating and solving integrated order batching and routing 
in multi-depot AGV-assisted mixed-shelves warehouses}




\author[mymainaddress]{Lin Xie\corref{mycorrespondingauthor}}
\cortext[mycorrespondingauthor]{Lin Xie}
\ead{xie@leuphana.de}
\ead[url]{http://www.leuphana.de/or/}
\author[mysecondaryaddress]{Hanyi Li}
\author[mymainaddress]{Laurin Luttmann}

\address[mymainaddress]{Leuphana University of L\"uneburg, Universit\"atallee 1, 21335 L\"uneburg}
\address[mysecondaryaddress]{Beijing Hanning Tech Co., Ltd}

\begin{abstract}
	Different retail and e-commerce companies are facing the challenge of assembling large numbers of time-critical picking orders that include both small-line and multi-line orders. To reduce unproductive picker working time as in traditional picker-to-parts warehousing systems, different solutions are proposed in the literature and in practice. For example, in a mixed-shelves storage policy, items of the same stock keeping unit are spread over several shelves in a warehouse; or automated guided vehicles (AGVs) are used to transport the picked items from the storage area to packing stations instead of human pickers. This is the first paper to combine both solutions, creating what we call \textit{AGV-assisted mixed-shelves picking systems}. We model the new integrated order batching and routing problem in such systems as an extended multi-depot vehicle routing problem with both three-index and two-commodity network flow formulations. Due to the complexity of the integrated problem, we develop a novel variable neighborhood search algorithm to solve the integrated problem more efficiently. We test our methods with different sizes of instances, and conclude that the mixed-shelves storage policy is more suitable than the usual storage policy in AGV-assisted mixed-shelves systems for orders with different sizes of order lines (saving up to 62\% on driving distances for AGVs). Our variable neighborhood search algorithm provides optimal solutions within an acceptable computational time.
\end{abstract}

\begin{keyword}
Logistics \sep  Order batching \sep Routing \sep Mixed-shelves storage \sep AGV-assisted picking
\end{keyword}

\end{frontmatter}


\section{Introduction}\label{sec:intro}
The most important and time-consuming task in a warehouse is the collection of items from their storage locations to fulfill customer orders. This process is called \textit{order picking}, which may constitute about 50--65\% of operating costs. Therefore order picking is considered the highest-priority area for productivity improvements (see \cite{de2007design}). 
In a traditional manual order picking system (also called a \textit{picker-to-parts system}), the pickers spend 50\% of their working time on the task of walking (see \cite{Tompkins.2010}; for an overview of manual order picking systems see \cite{de2007design}). The unproductive working times require the picker-to-parts system to have a large workforce, especially for companies which have millions of small-sized items in large warehouses, such as the retailers Amazon, Alibaba, Zara, Zalando and Walmart. Many of them provide both brick-and-mortar stores and online shops to create a seamless shopping experience for customers (omnichannel flexibility). Due to the diversity of online shops, we concentrate on small-sized orders with different order lines. Especially during the COVID-19 pandemic, online grocery sales are growing threefold faster (see \cite{fabric.2020}). There are increasing demands for alternative warehousing systems to increase the efficiency of order picking, for example, robot-based compact storage and retrieval systems and robotic mobile fulfillment systems (see more details in \cite{azadeh2019robotized}). Here we consider a relatively new warehousing concept that does not use expensive fixed hardware and can be easily and quickly implemented, called \textit{AGV-assisted picking} (see \cite{boysen2019warehousing}, \cite{azadeh2019robotized}). 


\subsection{AGV-assisted picking systems}
\begin{figure}[H]
	\centering
	\includegraphics[width=0.26\textwidth]{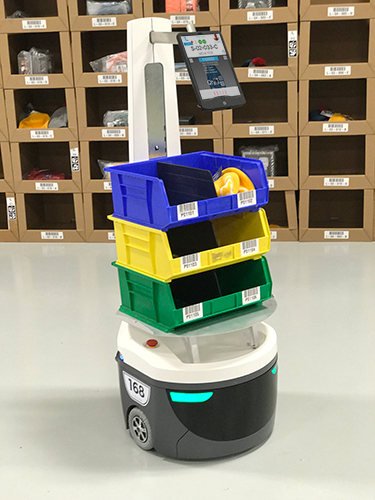}
	\includegraphics[width=0.3\textwidth]{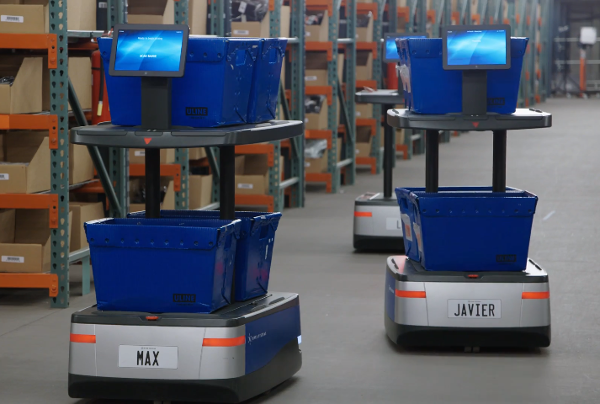}
	\includegraphics[width=0.3\textwidth]{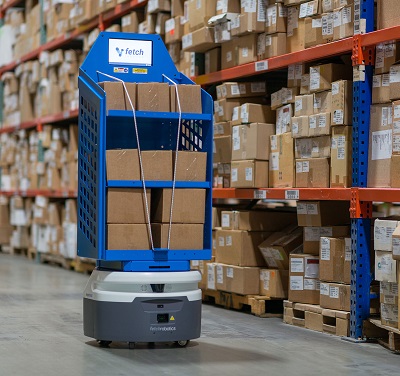}
	\caption{Companies who own an AGV-assisted picking system (from left to right: Locus Robotics, 6 River Systems, Fetch Robotics).}
	\label{fig:rafs_examples}
\end{figure}

As described in \cite{azadeh2019robotized}, there are different variants of this type of system. In this paper, we concentrate on the one that has been marketed to the warehousing industry, such as Locus Robotics, 6 River Systems and Fetch Robotics (see Figure~\ref{fig:rafs_examples}). 

\begin{figure}
	\centering
	\includegraphics[width=\textwidth]{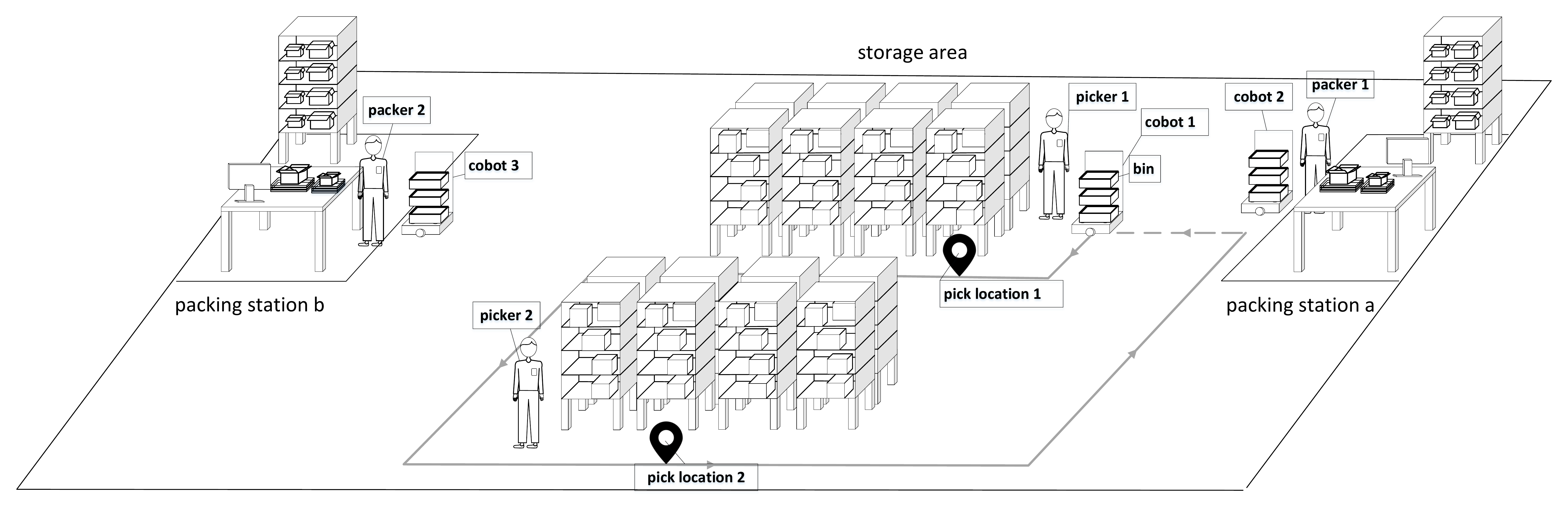}
	\caption{An example of a AGV-assisted picking system (with two packing stations).} 
	\label{fig:picking_refill_process}
\end{figure} 

In such systems, each human \textit{picker} works at his/her working area in the \textit{storage area} (containing shelves with items). For example, in Figure~\ref{fig:picking_refill_process}, picker 1 works in the upper area of the storage area while picker 2 works in the lower area. Each mobile robot (called a \textit{cobot}) leaves its packing station (also called a \textit{depot}), goes automatically to a pick location in the storage area and waits for a human picker to arrive. The picker identifies and picks up the ordered items from shelves and puts them into bins carried by the cobot (one \textit{bin} is used to temporarily store an order); after that, the cobot goes to its next picking location. The load of each cobot is limited. After collecting, each cobot goes back to its packing station to drop off the loaded bins, and the human \textit{packer} at that station packs the collected items. Each cobot has a dedicated packing station and it is charged there. In the example in Figure~\ref{fig:picking_refill_process}, cobot 1 leaves packing station $a$ and goes firstly to pick location 1 to wait for picker 1 to load items; after that, it goes to pick location 2 to wait for picker 2 to load items and returns to packing station $a$. If there is another cobot in front of station $a$, then it should wait in a \textit{queue} in front of the station. The explanation of all components in such a picking system can be found in Table~\ref{tab:components}. In different real-world warehouses, these depots might vary slightly as shown in Figure~\ref{fig:picking_refill_process}. For example, the depots might be drop-off locations, where the bins are put onto a conveyor system to transport them to packing stations. But each cobot still has a dedicated depot, which is its charging location. In the rest of the paper, the terms \textit{depot} and \textit{packing station} are synonyms. This system is called a \textit{follow-pick} system in \cite{lee2019robotics}, and \cite{boysen2019warehousing} call this model of operation the \textit{free-floating} policy. Although AGV-based picking systems are often applied for heavy and bulky items, there are also applications for small-sized items (see \cite{boysen2019warehousing}). The latter applications are also a focus of this paper.

In practice, AGV-assisted picking systems are becoming increasingly popular, since adding AGVs to an existing picker-to-parts system barely alters the basic fulfillment processes (easy and quick implementation). More than 80\% of all warehouses in Western Europe still follow the traditional picker-to-parts setup (see
\cite{de2012determining}). The other automated systems require, for example, expensive fixed hardware, such as automated storage systems in robot-based compact storage and retrieval systems. An AGV-assisted system is also easily scalable by adding or reducing pickers and AGVs to adapt to varying workloads, for example due to seasonal changes in customers' orders, such as end-of-season sales.

As described in \cite{azadeh2019robotized} and \cite{boysen2019warehousing}, AGV-assisted picking is a relatively new concept. The only publication about the free-floating system is in \cite{lee2019robotics}. They modeled the routing problem to support the cooperation of picker robots with transport robots and compared this with the cooperation of human pickers with transport robots (follow-pick). They concluded that human pickers are more suitable for identifying and picking up items. They also analyzed the relationship between layout configuration and robot system parameters. In contrast to that paper, we concentrate on integrated batching and routing in a mixed-shelves follow-pick system.

\subsection{Mixed-shelves picking systems}

As described in \cite{boysen2019warehousing} and \cite{weidinger2019picker}, mixed-shelves storage has been applied by many large-sized facilities in e-commerce companies, for example in the distribution centers of Amazon Europe and fashion retailer Zalando. In such storage, items of the same stock keeping unit (SKU) are spread over the shelves, so there are always some items of that SKU close by for picking. The main advantage of such storage is that it reduces unproductive working time. That is supported by \cite{weidinger2019picker} for traditional manual systems and \cite{Boysen.2017} for robotic mobile fulfillment systems. However, the literature shows that this storage is more suitable for online companies, which usually have many small-line orders (or small order-line demands). For cases with a large demand per SKU or multiple order lines per order, several storage locations might need to be visited by pickers to finish one order.

We apply mixed-shelves storage in an AGV-assisted picking system (called a \textit{multi-depot AGV-assisted mixed-shelves system}) for picking both single-line and multi-line small-sized orders and compare this with the reference storage policy, called \textit{dedicated storage}. In the dedicated storage policy, items of the same SKU are stored in only one shelf. 

\subsection{Operational decision problems}
The problems that occur in a multi-depot AGV-assisted mixed-shelves system are mostly similar to the problems in traditional manual order picking systems, including strategic (such as layout planning), tactical (such as zoning and storage assignment) and operational problems (batching, routing and job assignment) (see overview in \cite{van2018designing} for the traditional manual order picking systems). In addition, we have several depots in our system, and each cobot has its dedicated depot, i.e. it leaves from its depot and goes back to that depot. Therefore, we also need to decide in which depot an order will be handled to minimize driving distances while ensuring balanced batches among all depots (this will be discussed further in Section~\ref{sec:assump}).

We consider in this paper two integrated operational problems in a multi-depot AGV-assisted mixed-shelves system: customer orders are combined in several pick rounds (batches) in \textit{order batching} by deciding from which depot a batch is leaving, while the sequence of storage locations that cobots should visit to load all retrieved items from human pickers in each batch is determined in \textit{routing}. 

\subsection{Literature review}
In this subsection, the literature of picker routing in manual mixed-shelves warehouse and integrated order batching and routing will be described.

\subsection*{Picker routing in mixed-shelves picking} As shown in \cite{weidinger2019picker}, picker routing in a mixed-shelves warehouse is much more complex than in traditional environments: Multiple orders are concurrently picked by each picker (in our case: cobot), many alternative depots (in our case: packing stations) are available, and items of the same SKU are available in multiple shelves. \cite{daniels1998model} was the first publication to integrate the selection of alternative storage positions into picker routing. The same problem was also addressed in \cite{weidinger2018picker} for rectangular warehouses, and integrating alternative pick locations was proven to make this problem strongly NP-hard. 

The mixed-shelves storage in \cite{weidinger2019picker} is similar to the environment we consider, including:
\begin{itemize}
	\item Items of the same SKU are available in multiple shelves.
	\item The cobot / picking cart is able to carry multiple bins for different orders concurrently.
\end{itemize}
However, we do not have the same meaning of \textit{depots} as in \cite{weidinger2019picker}. They define depots as access points to the conveyor system, not the leaving and return locations of cobots. So in our case we have dedicated assignment of cobots to depots. Each cobot begins from its depot and ends at the same depot. This property makes our formulation similar to the vehicle routing problem with multiple depots. An overview of multi-depot vehicle routing problems can be found in \cite{montoya2015literature}. It is known that multi-depot vehicle routing problems are proven to be NP-hard (\cite{garey1979computers}). 

And we consider multiple pickers (in our case: cobots), rather than one picker as in \cite{weidinger2019picker}.
The following problems are mentioned in \cite{weidinger2019picker} for  considering multiple pickers.
\begin{itemize}
	\item They block each other in front of shelves (in this paper we assume that the aisles are wide enough for at least two (even three) cobots operating side by side (see \cite{lee2019robotics} for the same assumption)). 
	\item They block each other in front of depots (there is a queue in front of each depot; furthermore, due to the short drop-off time by removing bins from cobots, we assume that the blocking is minimal).
	\item They influence each other's inventory levels (this will be formulated in our model, and that might increase the complexity of the model).
\end{itemize}

\subsection*{Integrated batching and routing}
Due to the similar time horizon of the batching and routing decisions and the strong relationship between them, the integration of them is extensively studied in the literature; see an overview in \cite{van2019formulating}. The efficiency of the integration compared with the sequential approach to both of these problems is already shown in previous papers, such as in \cite{van2019formulating}. Therefore, we do not focus on such comparison in this paper.  

The differences in our problem from the existing integrated batching and routing problems for the traditional manual picking systems are summarized as follows. 
\begin{itemize}
	\item Due to mixed-shelves storage, the items of an SKU are available in several locations (shelves) (mixed storage, see third column of Table~\ref{tab:OverviewIntegratedBatchingAndRouting}). Therefore, the picking locations of an order should be determined during routing. 
	\item Multiple depots (packing stations) (see fourth column of Table~\ref{tab:OverviewIntegratedBatchingAndRouting}).
\end{itemize}

In the first integrated model from \cite{won2005joint}, they formulate it as the integrated bin packing and traveling salesman problem. The three-index formulations with sub-tour elimination constraints are widely used in the literature (see Table~\ref{tab:OverviewIntegratedBatchingAndRouting}). It is worth mentioning that \cite{chen2015efficient}, \cite{scholz2017order} and \cite{van2019formulating} also use three-index formulation to integrate batching and routing with another operational problem.
\begin{table}[H]
	\centering
	\begin{adjustbox}{width=\columnwidth,center}
	\begin{tabular}{ l l l l l l }
		\hline
		& Formu- & Mixed & Multiple  &Min. & Heuristic\\
		& lation& storage & depots & dist. &\\
		\hline
		\cite{won2005joint} & 1 &  &   & & FCFS+2-Opt\\
		\cite{tsai2008using} & 1 &  &  &  X& GA\\
		\cite{ene2012storage} & 1 &  &  &  & GA\\
		\cite{kulak2012joint} & 1 &  &  &  X& Cluster-based TS\\
		\cite{matusiak2014fast} & 1 & &   & X& SA+A*\\
		\cite{cheng2015using} & 1 &  & &   X&PSO+ACO\\
		\cite{lin2016joint} & 1 &  &  &  X& PSO\\
		\cite{li2017joint} & 1 &  &  &  X&Similarity+ACO\\
		\cite{valle2017optimally} & 1 &  &  &  X&-\\
		\cite{briant2020efficient} & & &  & & CG\\
		\textbf{this paper} & 1, 2 & X &X &  X & VNS\\
		\hline
	\end{tabular}
\end{adjustbox}
	\captionof{table}{The literature about integrated batching and routing. 1: three-index formulation; 2: two-commodity network flow formulation; FCFS: first come, first serve; GA: genetic algorithm; TS: tabu search; SA: simulated annealing; PSO: particle swarm optimization; ACO: ant colony optimization; CG: column generation; VNS: variable neighborhood search.}
	\label{tab:OverviewIntegratedBatchingAndRouting}
\end{table}

The two-commodity network flow formulation for solving the traveling salesman problem was first introduced by \cite{finke1984two}. This formulation was proposed by \cite{baldacci2004exact} to solve the capacitated vehicle routing problem and was extended by \cite{ramos2019multi} to solve the multi-depot vehicle routing problem, and it was shown that the flow formulation provides better performance to achieve lower gaps within a smaller computational time. In our work, we extend the formulation in \cite{ramos2019multi} for integrating vehicle routing with order batching in a mixed-shelves warehouse. In other words, we allow each location to be visited more than once and the demand in each location is unknown (but the maximum is given), and we enable batching of orders. More about this formulation can be found in Subsection~\ref{subsec:two-commodity}. 

Due to the complexity of the integration, the problem is solved in the literature mostly using metaheuristics (see sixth column of Table~\ref{tab:OverviewIntegratedBatchingAndRouting}). There are two main reasons why we apply variable neighborhood search to this integrated problem. First, \cite{li2017joint} mentioned that there are several problems with existing applied metaheuristics in the literature for the integrated batching and routing problem, such as the bad performance caused by not using the full capacity of the vehicle in the batching process, or the long computational time. It is not clear which of the heuristics proposed so far in the literature is the best one. Second, our mixed-shelves storage increases the complexity of the routing problem, and variable neighborhood search has been successfully applied to solve vehicle routing problems, such as large-scale capacitated vehicle routing in \cite{kytojoki2007efficient}.



\subsection{Contributions and paper structure}
Based on the literature review in the previous subsection, we summarize the main contributions of this paper as follows:
\begin{itemize}
	\item We formulate the first integrated batching and routing problem in a multi-depot AGV-assisted mixed-shelves warehouse with a new two-commodity network flow model and compare it with the extended classical three-index model. 
	\item Due to the high complexity of the integration, a novel variable neighborhood search is developed to provide good solutions within shorter computational times. 
\end{itemize}

The remainder of this paper is organized as follows. Section~\ref{sec:assump} describes the assumptions we used, while Sections~\ref{sec:model} and \ref{sec:heuristic} describe the different mathematical formulations and implemented metaheuristic for our integrated order batching and robot routing respectively. Computational results are described in Section~\ref{sec:results}. Finally, Section~\ref{sec:conclusions} presents conclusions and opportunities for future research.
 



\section{Assumption for the mathematical models} \label{sec:assump}
The following assumptions are made for the formulation of the integrated problem in Section~\ref{sec:model}.

\paragraph{Objective} We use the objective of minimizing distances for the reason that this objective allows us to simplify our integrated problem in such a way that the waiting time of a cobot for a human picker to arrive at the given picking location can be ignored. Therefore, the cobots act like classical pickers in a warehouse, i.e. they move from their depots to different picking locations to collect items for orders. In our model, we don't analyze the cooperation between human pickers and cobots at picking locations. This objective is also widely used in the literature for the integrated batching and routing problem in traditional manual systems (see fifth column of Table~\ref{tab:OverviewIntegratedBatchingAndRouting}). It is worth mentioning that this assumes that minimizing distances has the same meaning as minimizing total picking time in the literature (\cite{li2017joint}, \cite{tsai2008using} and \cite{won2005joint}). This is supported by the assertion in \cite{de2007design} that the travel time is an increasing function of the travel distance, so that minimizing the latter is considered a primary objective in warehouse design and optimization. However, this is no longer true in our problem, if we consider the waiting time of a cobot for a human picker in the mixed-shelves storage. The waiting time of a cobot for a human picker in the storage area can include the travel time of the human picker to the pick location and the waiting time in the queue (if other cobots arrive earlier than the current one and the picking has not proceeded, they automatically form a queue at the pick location). A pick location with a shorter driving distance does not mean that the picking time at that pick location is shorter in our problem. Another pick location with items of the same SKU with a longer driving distance might be preferred if the waiting time of the cobot is shorter. 
\paragraph{Pickers in the storage area} We don't analyze the number of pickers and size of their corresponding working areas (zones). 
\paragraph{A given set of orders} The small-line and multi-line orders of customers are given as input for the models. No new incoming orders are considered.
\paragraph{Demand nodes} In order to formulate the problem as a vehicle routing problem, we need to define a set of demand nodes. So it is possible to define each item as a node. But there are a large number of nodes, and all items on one shelf share the same distances to items on another shelf. Therefore we use a shelf as an aggregated node for all items on it. 
\paragraph{Picking list of an order} Each customer's order includes several SKUs with their given ordered quantities. Based on the property of mixed-shelves storage, items of the same SKU might be picked from different shelves. To simplify the modeling, we divide customer orders into a picking list. For example, we have order 1, including two units of SKU $a$, one unit of SKU $b$ and three units of SKU $c$:  $o_1={(a,2),(b,1),(c,3)}$ (SKU, ordered quantity); then we have a picking list for $o_1$: ${a_1,a_2,b_1,c_1,c_2,c_3}$ in our models. In this way, two items of SKU $a$ can be picked from different shelves.
\paragraph{No charging} We assume that all cobots we send out have enough remaining percentage of battery life for one tour. When a cobot visits its depot, it will be re-charged at a rate per unit of time. 
\paragraph{Aisle congestion} We assume that the aisles are large enough for two or even three cobots moving side by side (the same assumption is common in the literature, such as in \cite{lee2019robotics}), so the aisle congestion is minimal, since one cobot can go around the boundaries (such as: picker, another cobot).
\paragraph{Load capacity of batches} In integrated batching and routing, we can determine a tour for a cobot for each batch. So the load capacity of each batch is equal to the load capacity of a cobot. We assume that the payloads of all cobots are the same. 
\paragraph{Same number of cobots in each depot} We assume that the same number of cobots are available in each depot.
\paragraph{A bin/sub-bin for each order} Each cobot can carry several bins, and a bin can be separated into several small sub-bins for small-sized items (see, for example, Locus Robotics in Figure~\ref{fig:rafs_examples}). A human operator in a depot can rearrange the bins by removing or adding some divisions before the cobot leaves the depot. The rearrangement depends on the sizes of the orders to be transported in the next batch by this cobot. Due to the possibility of this rearrangement, we assume that the batching is limited by the load capacity of a cobot, not by the number of orders/bins as in \cite{valle2017optimally}.
\paragraph{Balanced batches of different depots} 
We have multiple depots in this work. Due to our objective of minimizing distances, it is possible that we have extremely unbalanced numbers of batches among depots. This is caused by the placement of depots and the shelves that include ordered items. For example, where ten batches are assigned to depot $a$ while only one batch is assigned to depot $b$, since most of the required shelves are directly in front of depot $a$. This assignment reduces the distances; however, it causes unbalanced batches among stations and consequently a bad makespan for a given set of orders. The main reason is that the more batches are assigned to one depot compared to others, the more waiting time of the human packer is needed at that depot for additional batches (we assume the same number of cobots for each depot), and also the more preparation time (e.g. labeling) is needed for each additional batch. Furthermore, we don't use the balanced distribution of the number of orders among different depots to achieve good makespan, since we have small-line and multi-line orders. For example, we have two depots ($a$ and $b$; each owns one cobot), six orders with respective weights 4, 1, 1, 1, 1 and 1 (the first one with four items while the others with one item); the capacity of each batch is limited to 5 of the total weight. We assign the first three orders to depot $a$ and the later three ones to depot $b$ to achieve balanced distribution of the number of orders among two depots. Two batches are needed leaving from depot $a$, and one is needed leaving from depot $b$. The human packer at depot $a$ needs to wait the second batch. If we use balanced batches instead, we will assign the first order to depot $a$ and all others to depot $b$, and we need only one batch leaving from each depot. As shown in this example, if we have items with similar weights, we can get a similar number of items within each batch (5 items at depot $a$ vs. 4 items at depot $b$) if they are as close to capacity as possible.
\paragraph{No refill in parallel to picking} We assume that our optimization begins after refill operations and the inventory has enough items for picking.


\section{Mathematical model for integrated order batching and cobot routing 
} \label{sec:model}
Our integrated order batching and cobot routing can be reduced to the capacitated vehicle routing problem if there is only one depot, with each shelf including only items of one SKU, single-line orders (each with a distinct SKU) and a given number of batches. We assume each shelf corresponds to a customer and each batch corresponds to a vehicle, so each customer is assigned to exactly one vehicle. As the capacitated vehicle routing problem is known to be NP-hard (see, for example, \cite{baldacci2007recent}), so is our problem.


In this section, we firstly model the integrated order batching and routing problem as a three-index formulation in Subsection~\ref{subsec:three-index}, and then as a two-commodity network flow formulation in Subsection~\ref{subsec:two-commodity}.
\newcommand{\SetOfShevlesIncludePickItem}[1]{\mathcal{V}_{#1}^\text{S}}
\newcommand{\SetOfShevlesNeedRefillItem}[1]{\mathcal{V}_{#1}^\text{SF}}
\newcommand{\SetOfShelfNodes}{\mathcal{V}^\text{S}}
\newcommand{\SetOfStationNodes}{\mathcal{V}^\text{D}}
\newcommand{\SetOfStationNodesCopy}{\mathcal{V}^\text{D'}}
\newcommand{\SetOfRobotsByStation}[1]{\mathcal{R}_{#1}}
\newcommand{\SetOfBatchesForPicking}{\mathcal{B}}
\newcommand{\SetOfPickingRobotsByStation}[1]{\mathcal{R}_{#1}}
\newcommand{\SetOfPickingBatchesByStation}[1]{\mathcal{B}_{#1}}
\newcommand{\SetOfPickingRobots}{\mathcal{R}}
\newcommand{\SetOfRefillRobots}{\mathcal{R}^\text{F}}
\newcommand{\SetOfPickingItems}{\mathcal{P}}
\newcommand{\SetOfRefillItems}{\mathcal{P}^\text{F}}
\newcommand{\SetOfItemsInShelf}[1]{\mathcal{P}_{#1}}
\newcommand{\SetOfPickingItemsInShelf}[1]{\mathcal{P}_{#1}}
\newcommand{\SetOfRefillItemsInShelf}[1]{\mathcal{P}_{#1}^\text{F}}
\newcommand{\SetOfOrders}{\mathcal{O}}
\newcommand{\SetOfItemsInOrder}[1]{\mathcal{P}_{#1}}
\newcommand{\batchsize}[1]{\mathcal{B}_{#1}}
\newcommand{\OriginalDepotOfCopy}[1]{m_{#1}}

\subsection{Three-index formulation} \label{subsec:three-index}
The following model is based on the integrated batching and routing in \cite{cheng2015using}. In contrast to that, we consider mixed-shelves storage and multiple depots. We need one more index for depots in our model so that we can know from which depot a batch begins. Note that the batches leaving from a depot obtained from our model are independent of each other, i.e. the sequence of batches is irrelevant. So, a set of cobots can be sent out to collect items.


The sets, parameters and variables are as follows: 

\hspace{-0.3cm}\textbf{Sets}:\\
\begin{tabularx}{\textwidth}{>{\hsize=0.45\hsize}X >{\hsize=1.55\hsize}X}
	$\mathcal{V}$ & Set of nodes, including shelves $\SetOfShelfNodes$ (including at least one ordered item) and depots
	 $\SetOfStationNodes$ ($\mathcal{V}=\SetOfShelfNodes \cup \SetOfStationNodes$)\\
	$\mathcal{E}$ & Set of edges $(i,j) \in \mathcal{E} \ \ \forall i,j \in \mathcal{V}$\\
	$\SetOfPickingItems$ & Set of physical items for picking\\
	$\SetOfShevlesIncludePickItem{p}$ & Set of shelves including item $p \in \SetOfPickingItems$ \\
	$\SetOfOrders$ & Set of orders\\
	$\SetOfItemsInOrder{o}$ & Set of physical items in the picking list for order $o$\\
	$\batchsize{d}$ & Set of batches leaving from depot $d \in \SetOfStationNodes$ \\	
\end{tabularx}
\textbf{Parameters}:\\
\begin{tabularx}{\textwidth}{>{\hsize=0.45\hsize}X >{\hsize=1.55\hsize}X}
	$d_{ij}$ & Distance between two nodes $i,j \in \mathcal{V}$\\
	$c$ & Maximum payload of a cobot\\
	$w_p$ & Weight of item $p \in \SetOfPickingItems$\\
	$n_{ps}$ & Number of items sharing the same SKU as the physical item $p$ at shelf node $s \in \SetOfShelfNodes$\\
\end{tabularx}
\textbf{Decision variables}:\\
\begin{tabularx}{\textwidth}{>{\hsize=0.2\hsize}X >{\hsize=1.8\hsize}X}
	$x_{ijbd}$ & 
	$
	\left\{
	\begin{array}{cl}
	1 , & \mbox{Node $j \in \mathcal{V}$ has been visited after node $i \in \mathcal{V}$ in batch $b \in \batchsize{d}$ leaving from depot $d$}\\
	0 , & \mbox{else}
	\end{array}
	\right.
	$
	\\
	$z_{psbd}$ & 
	$
	\left\{
	\begin{array}{cl}
	1 , & \mbox{Item $p$ has been collected at node $s \in \SetOfShevlesIncludePickItem{p}$ in batch $b \in \batchsize{d}$ leaving from depot $d$}\\
	0 , & \mbox{else}
	\end{array}
	\right.
	$
	\\
	$\omega_{obd}$ & 
	$
	\left\{
	\begin{array}{cl}
	1 , & \mbox{Order $o \in \SetOfOrders$ is collected in batch $b \in \batchsize{d}$ leaving from depot $d$}\\
	0 , & \mbox{else}
	\end{array}
	\right.
	$
	\\
\end{tabularx}
\\

As explained in Section~\ref{sec:assump}, we want to have a balanced number of batches among depots. Therefore, we can calculate the smallest upper bound on the number of batches $|\batchsize{d}|$ for all depots $d$ as follows: 
\begin{equation}
|\batchsize{d}|=\left\lceil \frac{\sum_{p \in \SetOfPickingItems} w_p}{|\SetOfStationNodes|*c} \right\rceil \label{eq:defBr}
\end{equation}
It is possible that a batch leaving from one depot will not be used according to our minimization objective function. It is also possible that there will be no feasible solutions to achieve $|\batchsize{d}|$. In this case, the number of batches will be increased by one or two, etc until we get a feasible solution.
\begin{description}
	\item[Goal:] \textbf{Minimizing traveled distances of all batches}
\end{description}
\begin{align}
F(3.1) \ \ \ 
z(F(3.1)) &= \text{Min} \quad  \sum_{d \in \SetOfStationNodes} \sum_{b \in \batchsize{d}} \sum_{(i,j) \in \mathcal{E}} d_{ij} * x_{ijbd}
\label{obj:distances} \\
\text{s.t.} \quad  \sum\limits_{j \in \mathcal{V}} x_{ijbd} &=  \sum\limits_{j \in \mathcal{V}} x_{jibd} \leq 1, \; \forall i \in \mathcal{V}, d \in \SetOfStationNodes, b \in \batchsize{d} \label{eq1:flowconservation}  \\
 \sum\limits_{p \in \SetOfPickingItems} \sum\limits_{s \in \SetOfShevlesIncludePickItem{p}} z_{psbd} * w_p &\leq c, \; \forall d \in \SetOfStationNodes, b \in \batchsize{d}\label{eq7:pickingcap}\\
\sum\limits_{i \in S} \sum\limits_{j \in S} x_{ijbd} &\leq \lvert S\rvert - 1, \; \forall S \subset \SetOfShelfNodes, \lvert S\rvert \geq 0, \forall d \in \SetOfStationNodes, b \in \batchsize{d} \label{eq9:sec}\\
 \sum\limits_{s \in \SetOfShevlesIncludePickItem{p}} \sum\limits_{d \in \SetOfStationNodes} \sum\limits_{b \in \batchsize{d}} z_{psbd}&=1, \; \forall p \in \SetOfPickingItems \label{eq3:exactonepicking}\\
 \sum\limits_{p \in \SetOfItemsInOrder{o}} \sum\limits_{s \in \SetOfShevlesIncludePickItem{p}} z_{psbd}& =\omega_{obd}*|\SetOfItemsInOrder{o}|, \; \forall o \in \SetOfOrders, d \in \SetOfStationNodes, b \in \batchsize{d} \label{eq:nonsplitting}\\
 z_{psbd} &\leq \sum\limits_{j \in \mathcal{V}} x_{sjbd}, \; \forall p \in \SetOfPickingItems, s \in \SetOfShevlesIncludePickItem{p}, d \in \SetOfStationNodes, b \in \batchsize{d} \label{eq5:relationxzpicking}\\
 \sum\limits_{d \in \SetOfStationNodes, b \in \batchsize{d}} z_{psbd} &\leq n_{ps}, \; \forall p \in \SetOfPickingItems, s \in \SetOfShevlesIncludePickItem{p} \label{eq:limitednumberofitemspershelf}\\
  x_{ijbd} &\in \lbrace 0, 1 \rbrace, \; \forall (i,j) \in \mathcal{E}, d \in \SetOfStationNodes, b \in \batchsize{d} \label{eq11:defx}\\
 z_{psbd} &\in \lbrace 0, 1 \rbrace, \; \forall p \in \mathcal{P}, d \in \SetOfStationNodes, b \in \batchsize{d}, s \in \SetOfShevlesIncludePickItem{p} \label{eq12:defz}\\
\omega_{obd} &\in \lbrace 0, 1 \rbrace, \; \forall o \in \SetOfOrders, d \in \SetOfStationNodes, b \in \batchsize{d} \label{eq:definition_omega_ob}
\end{align}
\begin{description}
	\item[A.] \textbf{Same constraints as in classical multi-depot vehicle routing problems} 
\end{description}
\hspace{-0.01cm}Constraint set \eqref{eq1:flowconservation} indicates the flow conservation and ensures that each batch leaves its depot once at the most, and if it leaves its depot, it should come back to that depot at the end of the tour. 
Furthermore, the capacity of each batch is limited in constraint set \eqref{eq7:pickingcap}. Constraint set \eqref{eq9:sec} eliminates subtours. Similar constraints can be found in, for example, \cite{ramos2019multi} for multi-depot vehicle routing problems.

\begin{description}
	\item[B.] \textbf{Constraints about order batching}
\end{description}
\paragraph{Each item must be assigned to a batch} Each item in the picking list should be collected exactly once in one batch in constraint set \eqref{eq3:exactonepicking}.
\paragraph{No splitting}
 Constraint set \eqref{eq:nonsplitting} ensures that all items for an order are collected in precisely one batch assigned to a cobot.

\paragraph{Set of constraints due to mixed-shelves storage}
Constraint set \eqref{eq5:relationxzpicking} ensures that if one item $p$ in shelf $s$ is collected in batch $b$, shelf $s$ should be visited in that batch. Constraints \eqref{eq:limitednumberofitemspershelf} ensure that the available stock of items in a shelf is not exceeded in all batches.

\begin{description}
	\item[D.] \textbf{Definition of variables}
\end{description}
\hspace{-0.01cm}The binary definitions of variables are given in constraint sets \eqref{eq11:defx}--\eqref{eq:definition_omega_ob}.

\subsection{Two-commodity network flow formulation} \label{subsec:two-commodity}
The model about the two-commodity network flow formulation is based on the formulation presented by \cite{ramos2019multi} for solving the multi-depot vehicle routing problem. Constraint set \eqref{eq1:flowconservation} and the constraint sets about order batching (\eqref{eq3:exactonepicking}--\eqref{eq:definition_omega_ob}) are extended with a set of copy depots and new flow variables. We will show you an example to illustrate the idea of the flow formulation. First we need some additional definitions of sets, parameters and variables.\\
\textbf{Additional sets}:\\
\\
\begin{tabularx}{\textwidth}{>{\hsize=0.45\hsize}X >{\hsize=1.55\hsize}X}
	$\SetOfStationNodesCopy$ & Set of copy of depots ($\mathcal{V}'=\SetOfShelfNodes \cup \SetOfStationNodes \cup \SetOfStationNodesCopy$)\\
	$\SetOfItemsInShelf{s}$ & Set of all items in shelf $s \in \SetOfShelfNodes$\\
	$\mathcal{E}'$ & Set of edges, $\mathcal{E}'=\mathcal{E} \cup \{(i,j): i \in \SetOfShelfNodes, j \in \SetOfStationNodesCopy\}$
\end{tabularx}
\textbf{Additional parameter}:\\
\\
\begin{tabularx}{\textwidth}{>{\hsize=0.45\hsize}X >{\hsize=1.55\hsize}X}
	$\OriginalDepotOfCopy{i}$ & Real depot of the copy depot $i \in \SetOfStationNodesCopy$\\
\end{tabularx}
\textbf{Additional decision variables}:\\
\\
\begin{tabularx}{\textwidth}{>{\hsize=0.2\hsize}X >{\hsize=1.8\hsize}X}
	$y_{ijbd}$ & A flow variable representing the batch load when the path from $i$ to $j$ is travelled in the $b$-th batch leaving from depot $d$. The flow $y_{jibd}$ represents the current empty space in the $b$-th batch leaving from depot $d$.
	\\
\end{tabularx}
\begin{figure}[h]
	\centering
	\includegraphics[width=\textwidth]{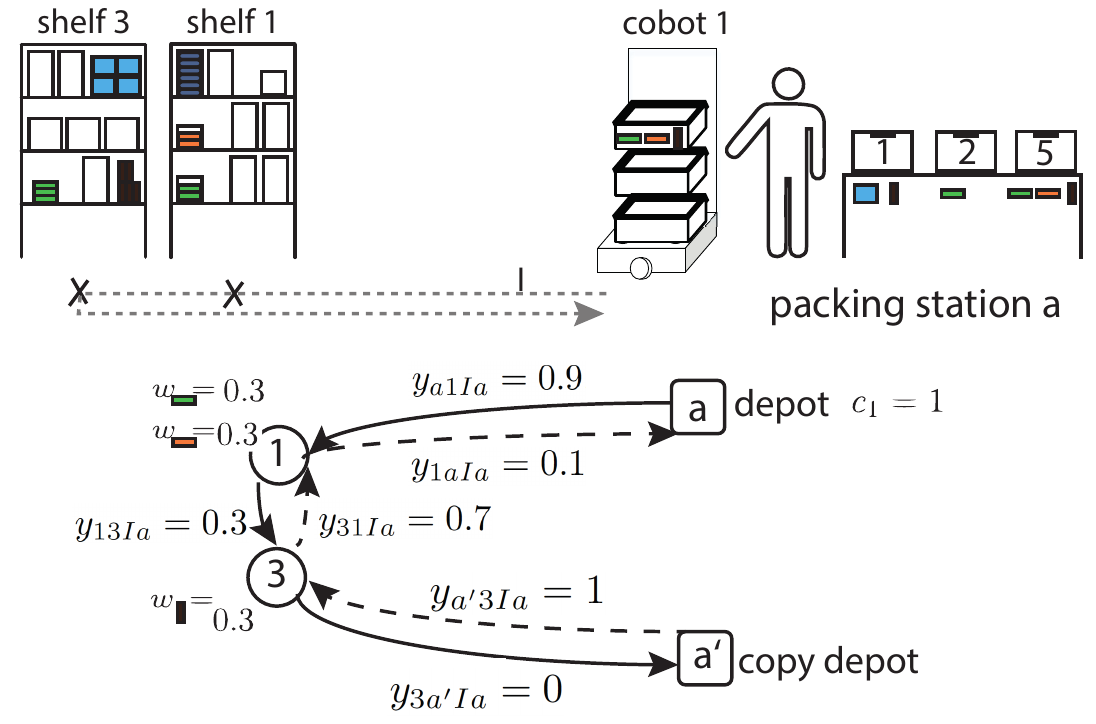}
	\caption{Flow paths for batch $I$ for picking three items.}
	\label{fig:example_two_commodity}
\end{figure} 
\paragraph{Example} Figure~\ref{fig:example_two_commodity} shows an example of batch $I$ leaving from depot $a$ limited by the batch capacity $c_I=1$. Three orders are assigned to that depot. They include different items, whose sizes are equal (0.3). The items sharing the same SKU are illustrated with the same color and shape. For example, the items indicated by a green rectangle are included in orders 2 and 5, and they are stored in shelves 1 and 3. In the flow formulation, the route for batch $I$ is defined by two paths: one from the real depot $a$ to the copy depot $a'$ ($a \rightarrow 1 \rightarrow 3 \rightarrow a'$), defined by $y_{ijrb}$ (representing the batch load); while the reserve path from the copy depot to the real depot ($a' \rightarrow 3 \rightarrow 1 \rightarrow a$) is defined by $y_{jibd}$ (representing the empty space in the batch). For example, the flow $y_{a1Ia}$ indicates that batch I leaves the depot with a load equal to the total weights of the three picking items, 0.9; while the flow $y_{a'3Ia}=1$ means that the batch arrives empty at the depot after dropping off items. For each edge $(i,j)$ of the route within a batch $b$ leaving from depot $d$, we have $y_{ijbd}+y_{jibd}=c$, for example, $y_{13Ia}+y_{31Ia}=0.3+0.7=1$.

\begin{description}
	\item[Goal:] \textbf{Minimizing traveled distances for all batches}
\end{description}
\begin{align}
(F(3.2)) \ \ \ z(F(3.2)) &=\text{Min} \quad  \sum_{d \in \SetOfStationNodes} \sum_{b \in \batchsize{d}} \sum_{(i,j) \in \mathcal{E}} d_{ij} * x_{ijbd}
\tag{\ref{obj:distances}} \nonumber\\
\text{s.t.} \quad 
&\eqref{eq1:flowconservation}, \eqref{eq3:exactonepicking}-\eqref{eq:definition_omega_ob} \nonumber \\
\sum_{j \in \mathcal{V}'} (y_{jsbd} - y_{sjbd}) &= 2 * \sum_{p \in \SetOfPickingItemsInShelf{s}} (w_p * z_{psbd}), \; \forall s \in \SetOfShelfNodes, d \in \SetOfStationNodes, b \in \batchsize{b}  \label{eq17:inflow-outflow}\\
\sum\limits_{d \in \SetOfStationNodes} \sum\limits_{j \in \SetOfShelfNodes} \sum\limits_{b \in \batchsize{d}} y_{djbd} &= \sum_{p \in \SetOfPickingItems} w_p \label{eq18:totaloutflow}\\
 \sum\limits_{d \in \SetOfStationNodes} \sum\limits_{j \in \SetOfShelfNodes} \sum\limits_{b \in \batchsize{d}} y_{jdbd} &\leq  c * \sum_{d \in \SetOfStationNodes} |\batchsize{d}| - \sum_{p \in \SetOfPickingItems} w_p \label{eq19:totalinflow_realdepot}\\
\sum_{j \in \SetOfShelfNodes} y_{ijbd} &\leq c, \; \forall i \in \SetOfStationNodesCopy, d=\OriginalDepotOfCopy{i}, b \in \batchsize{d} \label{eq20:outflow_copydepot}\\
y_{ijbd} + y_{jibd} &= c * x_{ijbd}, \; \forall i, j \in \mathcal{V}', d \in \SetOfStationNodes, b \in \batchsize{d} \label{eq23:capacity}\\
y_{sjbd} &\leq M * \sum_{p \in \SetOfItemsInShelf{s}} z_{psbd}, \; \forall s \in \SetOfShelfNodes, j \in \mathcal{V}', d \in \SetOfStationNodes, b \in \batchsize{d} \label{eq24:relationyz}\\
\sum_{j \in \SetOfShelfNodes} x_{ijbd} &= 0, \; \forall i \in \SetOfStationNodesCopy, d \neq \OriginalDepotOfCopy{i}, b \in \batchsize{d} \label{eq25:othercopydepot}\\
\sum_{j \in \SetOfShelfNodes} x_{jibd} &= 0, \; \forall i \in \SetOfStationNodes, d \neq \OriginalDepotOfCopy{i}, b \in \batchsize{d} \label{eq26:otherrealdepot}\\
y_{isbd} \geq \sum_{p \in \SetOfItemsInShelf{s}} w_p*z_{psbd} &- c*(1-x_{isbd}), \; \forall s \in \SetOfShelfNodes, i \in \mathcal{V}', d \in \SetOfStationNodes, b \in \batchsize{d} \label{eq:valid_inequality}\\
y_{ijbd} &\geq 0, \; \forall i, j \in \mathcal{V}', d \in \SetOfStationNodes, b \in \batchsize{d} \label{eq27:defy}
\end{align}

The objective function \eqref{obj:distances} is the same as in F(3.1), while some constraint sets are retained as in F(3.1), including \eqref{eq1:flowconservation} and \eqref{eq3:exactonepicking}--\eqref{eq:definition_omega_ob}. For each batch, the inflow minus outflow at each shelf node is equal to twice the weights of picked items on that shelf (see constraint set \eqref{eq17:inflow-outflow}). The total outflow from the real depots for all batches indicates the total load, which should be equal to the weight of all items (see constraint \eqref{eq18:totaloutflow}), while the total inflow to the real depots should be smaller than or equal to the residual capacity of the used batches (see constraint \eqref{eq19:totalinflow_realdepot}). 
Constraint set \eqref{eq20:outflow_copydepot} ensures that the outflow of each batch beginning with a copy depot $i$ is less than or equal to that batch's capacity. Constraint set \eqref{eq23:capacity} ensures that the inflow plus outflow of each node equals the capacity of the batch. Constraint set \eqref{eq24:relationyz} guarantees that the flow variable of batch $b$ from shelf $s$ is set to zero; if not, a single item $p$ from shelf $s$ is picked in batch $b$ leaving from depot $d$. Each flow variable is limited by the batch capacity, i.e. $M = c$. Finally, constraint sets \eqref{eq25:othercopydepot} and \eqref{eq26:otherrealdepot} jointly ensure that a batch cannot leave and return to a copy/real depot other than its home depot. Constraint set \eqref{eq:valid_inequality} defines valid inequalities to tighten the formulation by reducing the search space. These constraints ensure that the remaining batch capacity before reaching shelf $s$ is larger than or equal to the weights of all picked items at shelf $s$. Constraint set \eqref{eq27:defy} defines new linear flow variables.

\section{Variable neighborhood search algorithm for integrated order batching and cobot routing} \label{sec:heuristic}
As described before, integrated batching and routing is very complex, so it is solved in the literature mostly using metaheuristics to get a good solution within a reasonable computational time (see Table~\ref{tab:OverviewIntegratedBatchingAndRouting}, the last column). In this work, we use variable neighborhood search (VNS) to solve this integrated problem. Variable neighborhood search is a relatively new metaheuristic (proposed by \cite{mladenovic1997variable}), and its basic idea is \textit{a systematic change of neighborhood both within a descent phase to find a local optimum and in a perturbation phase to get out of the corresponding valley} (according to \cite{gendreau2010handbook}). 
\begin{algorithm}[h]
	\DontPrintSemicolon
	\textit{Initialization}. Select the set of neighborhood structures $N_\kappa (\kappa=1,...,\kappa_{max})$ that will used in the search; find an initial solution $x$; choose a stopping condition\;
	\textit{Repeat} the following until the stopping condition is met:
	\begin{enumerate}
		\item Set $\kappa \leftarrow 1$;
		\item Repeat the following steps until $\kappa=\kappa_{max}$:
		\begin{enumerate}
			\item[(a)] \textit{Shaking}. Generate a point $x'$ at random from ${\kappa}^{th}$ neighborhood of $x$ ($x'\in N_\kappa(x)$);
			\item[(b)] \textit{Local search}. Apply some local search method with $x'$ as initial solution; denote with $x''$ the so obtained local optimum;
			\item[(c)] \textit{Move or not}. If the local optimum $x''$ is better than the incumbent, move there $x\leftarrow x''$, and continue the search with $N_1 (\kappa \leftarrow 1)$; otherwise, set $\kappa \leftarrow \kappa+1$
		\end{enumerate}
	\end{enumerate}	
	\caption{Steps of the basic VNS (based on \cite{mladenovic1997variable})}
	\label{alg:vns}
\end{algorithm}

\newcommand{\SetOfOrdersForStation}[1]{\mathcal{O}_{#1}}
\newcommand{\SetOfOrdersCurrent}{\mathcal{O}_{cur}}
\newcommand{\SetOfItemsCurrent}{\mathcal{P}^\text{C}_{cur}}
\newcommand{\SetOfPickingItemsForStation}[1]{\mathcal{P}^\text{C}_{#1}}
\newcommand{\SetOfShevlesVisitedFromStation}[1]{\mathcal{V}_{#1}^\text{SC}}
\newcommand{\SetOfShevlesCurrent}{\mathcal{V}^\text{SC}_{cur}}
\newcommand{\SetOfShelvesSequenceForStation}[1]{\mathcal{V}_{#1}}

The steps of the basic VNS are shown in Algorithm~\ref{alg:vns}. Here, $N_\kappa (\kappa=1,...,\kappa_{max})$ is a set of neighborhood structures, which we will explain together with shaking in Subsection~\ref{subsec:shaking}. The stopping condition may be, for example, the number of iterations ($\gamma$) without improvement of $x^{\prime}$ or the maximum computational time $t_{max}$. It is worth mentioning that the algorithms we implemented in this section work for both mixed-shelves and dedicated storages.  Therefore, our algorithms can also be applied to the integrated batching and routing problem for the traditional manual picking systems.


In the following subsections, the implementation of each part of VNS for solving our integrated problem is described, including the building of an initial solution (Subsection~\ref{subsec:initial}), the shaking phase with the definition of the neighborhood structure (Subsection~\ref{subsec:shaking}) and the local search methods (Subsection~\ref{subsec:ls}).
\subsection{Initial solution}\label{subsec:initial}
\begin{algorithm}
	\SetAlgoLined
	\DontPrintSemicolon
	\KwResult{Initial solution $x$}
	\textbf{Function} \textsc{Initial}($\mathcal{O}, \, \mathcal{V}^D$)\;
	$x \gets \emptyset$\;
	$\mathcal{O}_d \gets \emptyset \quad \forall d \in \mathcal{V}^{\mathrm{D}}$\;
    $\hat{D}^o \gets \{\hat{d}_{od} \, | \,  d \in \mathcal{V}^{\mathrm{D}} \} \quad \forall \, o \in \mathcal{O} \,$\;
	$\mathcal{O}_d \gets$ \textsc{AssignOrdersToDepots}($\mathcal{O}, \, \mathcal{O}_d, \, \hat{D}^o) \quad \forall d \in \mathcal{V}^{\mathrm{D}}$\;
	\For{$ d \in  \mathcal{V}^D$}{
		$\mathcal{B}_d \gets$ \textsc{AssignOrdersToBatchesAndTours}($\SetOfOrdersForStation{d}$)\;
		$x \leftarrow x \cup \{\mathcal{B}_d\}$
	}
	\textbf{return} $x$\;
	\caption{Greedy algorithm to generate an initial solution}
	\label{alg:greedyInitial}
\end{algorithm}

In this subsection a greedy heuristic is presented in Algorithm~\ref{alg:greedyInitial}, which generates an initial solution $x$ for our integrated problem. 
First, the algorithm \textsc{AssignOrdersToDepots} determines the assignment of orders $o$ to depots $d$, while keeping the number of batches among depots balanced (see line 5 in Algorithm~\ref{alg:greedyInitial}). For this algorithm, a set of distance metrics $\hat{D}^o$ is needed, which we will describe together with algorithm \ref{alg:greedy_order_to_stations}. The set of orders assigned to the depot $d$ is denoted as $\mathcal{O}_d \subseteq \mathcal{O}$. Then, for each depot the orders assigned to it are bundled up into batches. This is done using the algorithm \textsc{AssignOrdersToBatchesAndTours} (line 7), which determines the set of batches $\mathcal{B}_d$ and the corresponding routes for a given depot $d$. The solution $x$ is a collection of all batches leaving from all depots and the corresponding routes for visiting shelves in batches.

\subsubsection{Assignment of orders to depots}
\begin{algorithm}[h]
	\SetAlgoLined
	\DontPrintSemicolon
	\KwResult{Set of orders $\mathcal{O}_d$ per station $d \in \mathcal{V}^{\mathrm{D}}$}
	\textbf{Function} \textsc{AssignOrdersToDepots}($\mathcal{O}, \, \mathcal{O}_d, \, \hat{D}^o$) \;
	$w_o := \sum_{p \in \mathcal{P}_o} w_p$ \tcp*[f]{calculate the weight of order $o$}\; 
	$w_d :=\sum_{o\in \mathcal{O}_d} w_o $ \;
    $\mu_o \gets \hat{D}^o_{(2)} - \hat{D}^o_{(1)} \quad \forall \, o \in \mathcal{O} $\;
    sort $\mathcal{O}$ according to $\mu_o$ in descending order \;
	\For{$o \in  \mathcal{O}$}{
        $d^* \gets \underset{d \in \mathcal{V}^D}{\mathrm{argmin}} (\hat{D}^o)$\;
		\uIf(\tcp*[f]{check balancing}){$\left \lceil \frac{ w_{d^*} + w_o }{c} \right \rceil \leq |\mathcal{B}_{d^*}|$}{
		    $\mathcal{O}_{d^*} \gets \mathcal{O}_{d^*} \cup o \, $; $ \mathcal{O} \gets \mathcal{O} \setminus o$\;
		}
		\Else{
			$\hat{D}^o \gets \hat{D}^o \setminus \hat{d}_{od^*}$\;
			$\mathcal{O}_d \gets \textsc{AssignOrdersToDepots}(\mathcal{O}, \, \mathcal{O}_d, \hat{D}^{o})$\;
		}
	}
	\textbf{return} $\mathcal{O}_d$\;
	\caption{Greedy assignment of orders to depots}
	\label{alg:greedy_order_to_stations}
\end{algorithm}

Algorithm~\ref{alg:greedy_order_to_stations} assigns orders to depots. For these types of assignment problems, it is common to first consider the distance from the object to be assigned (orders in our case) to all potential depots \citep{giosa2002new}. For our problem, this is already a challenge as each order consists of several items; further, each item can be stored in different shelves. Hence, the distance from an order to each depot is not well defined. To define a hypothetical loss for a given depot, we take for each item the distance of the shelf storing the item which is closest to the depot. These distances are then added up over all items of the order to get the hypothetical distance from that order to a depot. This definition of the distance from an order to a depot is formalized in Equation (\ref{eq:dist_order_station}):
\begin{align}
	\hat{d}_{od} = \sum_{p \in \mathcal{P}_o} \, \underset{j \in \mathcal{V}^{\mathrm{SC}}_p}{\mathrm{min}} \, (d_{ij}) \qquad \forall o \in \mathcal{O}, \, d \in \mathcal{V}^D,
	\label{eq:dist_order_station}
\end{align}
where we write $\hat{d}$ to indicate the hypothetical nature of the distance between order $o$ and depot $d$. These values are determined in line 4 of Algorithm \ref{alg:greedyInitial}. Moreover, we need to define the weights of orders as the sum of weights of all their items (line 2) together with the total weight of all orders assigned to each depot (line 3) in order to maintain a balance in the number of batches per depot.

Given the distances for all possible combinations of $(o,d)$, we follow the approach of \citet{giosa2002new} and calculate for each order the potential savings of assigning that order to its closest depot instead of the second closest (line 4 of algorithm \ref{alg:greedy_order_to_stations}). This is done in accordance with the following equation:

\begin{align}
	\mu_o = \hat{D}^o_{(2)} - \hat{D}^o_{(1)},
	\label{eq:calc_savings}
\end{align}
where $\hat{D}^o$ with subscripts $(1)$ and $(2)$ denotes the first- and second-order statistics respectively, i.e. $\hat{D}^o_{(2)}$ denotes the second smallest and $\hat{D}^o_{(1)}$ the smallest distance in $\hat{D}^o$.


The savings are then used to prioritize orders, i.e. to determine the sequence in which the orders are assigned to the depots (line 5). Due to the balancing constraint, orders cannot be assigned arbitrarily to stations, but must be assigned such that the number of batches at each station is as balanced as possible. Hence, orders where the second-best choice would be far worse than the best choice are prioritized and therefore assigned prior to orders where the difference is not too large. After choosing a depot $d^*$ for order $o$ in line 7, order $o$ can be assigned to a station in line 9 if the number of batches required to fulfill all assigned orders to this depot including $o$ is less than or equal to $|\mathcal{B}_{d^*}|$ as defined in Equation \ref{eq:defBr}. 

If an order cannot be assigned to the station closest to it due to the balancing constraint, the station is removed from the set of possible stations for that order (see line 11) and a new savings value for that order has to be determined. As the new savings value might affect the precedence of the remaining orders, this is done by recursively calling function \textsc{AssignOrdersToStations} with the reduced set of stations per order and only for the orders that are not assigned yet (line 12). This procedure is repeated until all orders are assigned. 

\subsubsection{Order batching and routing}
\begin{algorithm}
	\SetAlgoLined
	\DontPrintSemicolon
	\KwResult{Set of batches $\mathcal{B}_d$ with the size of $batchnr$}
	\textbf{Function} \textsc{AssignOrdersToBatchesAndTours}($\SetOfOrdersForStation{d}$)\;
	$batchnr \gets 1$\;
	\While{$\mathcal{O}_d \neq \emptyset$}{
		$weight \gets 0$; $tour \gets [d,d]$; $\mathcal{O}^{\mathrm{temp}}\leftarrow \emptyset$; $\mathcal{O}^{\mathrm{infeas}} \gets \emptyset$\;
		\While(\tcp*[f]{keep batch capacity}){$weight <c$}{
			$o \gets MinDist(\mathcal{O}_d \setminus \mathcal{O}^{\mathrm{infeas}}, \, \mathcal{O}^{\mathrm{temp}})$\;
			\uIf(\tcp*[f]{check batch capacity}){$weight + w_o < c$}{
				$weight \gets weight + w_o$\;
				$\mathcal{O}^{\mathrm{temp}}\leftarrow \mathcal{O}^{\mathrm{temp}} \cup o$\;
				$tour\leftarrow$\textsc{GreedyTour}({$tour, o$}) \tcp*{update the tour}
				$\mathcal{O}_d \leftarrow \mathcal{O}_d \setminus o$\;
			}\Else{
			    $\mathcal{O}^{\mathrm{infeas}} \gets \mathcal{O}^{\mathrm{infeas}} \cup o$
			}
		}
		$\SetOfBatchesForPicking_{d, batchnr}\leftarrow (\mathcal{O}^{\mathrm{temp}},tour)$; $batchnr \gets batchnr+1$\;
	}
	\textbf{return} $\SetOfBatchesForPicking_{d}$\; 
	\caption{Create a set of batches and tours for orders handled in depot $d$}
	\label{alg:greedy_order_batch}
\end{algorithm}
\begin{algorithm}
	\SetAlgoLined
	\DontPrintSemicolon
	\KwResult{Updated $tour$ by picking order $o$}
	\textbf{Function} \textsc{GreedyTour}($tour := [s_1,...,s_n,s_1], \, o $)\;
	$d^{add}_{best} \gets \infty$\;
	\For(\tcp*[f]{iterate over items of new order}){$p \in \mathcal{P}^{\mathrm{C}}_{o} $ }{
		\For{$s \in  V^{\mathrm{SC}}_p$}{
			\If(\tcp*[f]{check if shelf still has units}){$n_{ps} > 0$}{
				$j \gets \underset{j \in \{2,...,|tour|-1\}}{\mathrm{argmin}}\left (d_{s_{j-1}, s} + d_{s, s_{j}} -  d_{s_{j-1}, s_{j}} \right)$ \;
				$d^{add}_j \gets d_{s_{j-1}, s} + d_{s, s_{j}} -  d_{s_{j-1}, s_{j}}$\;
				\If{$d^{add}_j < d^{add}_{best}$}{
					$d^{add}_{best} \gets d^{add}_j$; $s^* \gets s$; $j^* \gets j$
				}
			}
		}
		$tour \gets [s_1,...,s_{j^*-1}, s^*, s_{j^*},...,s_n, s_1]$ \tcp*{insert node} 
		$n_{ps^*} \gets n_{ps^*} - 1$\tcp*{decrease supply by one}
	} 
	\textbf{return} $tour$ 
	\caption{Create a greedy batch tour for an order}
	\label{alg:greedy_tour_mixed}
\end{algorithm}

In Algorithm~\ref{alg:greedy_order_batch}, the orders $o \in \mathcal{O}_d$ of a depot $d \in \SetOfStationNodes$ are assigned to batches in compliance with the capacity restrictions while also forming tours. To initialize a new batch (see line 4), the initial weight, tour and set of orders in the batch, $\mathcal{O}^{\mathrm{temp}}$, are defined. Moreover, a set to store infeasible orders that do not fit in the batch, $\mathcal{O}^{\mathrm{infeas}}$, is initialized. In line 6, an order from the set of unassigned and feasible orders is selected by calculating the minimum possible distance from all items of $o$ to any item of the orders 
in the batch. Accordingly, the function $MinDist(\cdot)$ implements the following equation:
\begin{align}
    o = \underset{o\in \mathcal{O}_d \setminus \mathcal{O}^{\mathrm{infeas}}}{\mathrm{argmin}} \left (\sum_{p \in \mathcal{P}_o} \mathrm{min} (d_{st} \, | \, s \in \mathcal{V}^{\mathrm{S}}_{p}, \, t \in tour ) \right ) ,
\end{align}
where a random order is chosen if $\mathcal{O}^{\mathrm{temp}}$ is empty.
If the determined order still fits into the batch, it is added to it (lines 7--9); otherwise it is added to the set of infeasible orders $\mathcal{O}^{\mathrm{infeas}}$ (lines 12--13). 
If the order is added to the batch, its items need to be assigned to a shelf for pickup. Since an item can be stored in many shelves, the shelf that minimizes the distance to any of the shelves visited in the current batch is determined. 
Therefore, in line 10, each time a new order is assigned to a batch, the tour of this (partial) batch is updated using the cheapest insertion procedure (see e.g. \citet{archetti2003reoptimizing}), which we describe in algorithm \ref{alg:greedy_tour_mixed}. This algorithm determines, for each new item that has to be added to the current tour and for each shelf this item is stored in, the position in the tour where the distance added is minimized (see line 6). The additional distance created by adding a new node $s$ between nodes $s_{j-1}$ and $s_j$ can easily be calculated using the formula in line 7.
In the end, the item is picked up at the shelf $s^*$ which adds the least distance to the current tour (determined in line 9), and this shelf is inserted at the previously determined minimum distance position (in line 13).  

\subsection{Shaking} \label{subsec:shaking}
\begin{algorithm}[H]
	\SetAlgoLined
	\DontPrintSemicolon
	\SetKwFor{RepTimes}{repeat}{times}{end}
	\caption{Shaking}
	\KwResult{Solution $x^\prime$ from neighborhood $\mathcal{N}_\kappa (x)$ of current solution $x$}
	\label{alg:shake}
	\textbf{Function} \textsc{Shaking}($x, \kappa$)\;
	Determine total weight $w_b$ per batch $b \in \mathcal{B}_d$ in $x$; $\mathcal{O}^{\mathrm{temp}} \gets \emptyset$\;
	\RepTimes{$\kappa$}{
		$a \gets \textsc{EliminateBatch}(\mathcal{B}_d)$\tcp*{destroy a batch}
		$\mathcal{B}_d \gets \mathcal{B}_d \setminus a$; $\mathcal{O}^{\mathrm{temp}} \gets \mathcal{O}^{\mathrm{temp}} \cup \mathcal{O}_a$\;
	}
	\For(\tcp*[f]{repair orders}){$o \in \mathcal{O}^{\mathrm{temp}}$}{
		\uIf(\tcp*[f]{add orders to a new batch}){$\nexists \, b\in \mathcal{B}_d : \, w_b+w_{o} \leq c$}{
			$b_{|\mathcal{B}_d|+1} \leftarrow o; \, \mathcal{B}_d \gets \mathcal{B}_d \cup b_{|\mathcal{B}_d|+1}$; $b \gets b_{|\mathcal{B}_d|+1}$ \;
		}
		\Else(\tcp*[f]{add orders to a existing batch}){
			$b \gets$ select a random batch $b \in \mathcal{B}_d$ for which $w_b+w_{o} \leq c$\;
			$w_b \gets w_b + w_{o}; \; \mathcal{O}_b \gets \mathcal{O}_b \cup o$\;
		} 
		$tour_b \gets \textsc{GreedyTour}(tour_b, o)$\;
	}
	\SetArgSty{textnormal}
	\If{\textbf{any} $|\mathcal{B}_d|>|\mathcal{B}_d|^* \quad \forall d \in \mathcal{V}^{\mathrm{D}}$}{
		$\mathcal{B}_d \gets \textsc{ReduceBatches}(\mathcal{B}_d)$\;
	}
\textbf{return} $x$ 
\end{algorithm}
After an initial solution is obtained using the algorithms described in the previous subsection, its neighborhood can be searched for better solutions. The \textsc{Shaking} operator is mainly used to search increasingly distant neighborhoods of the current incumbent solution, which are then searched for a local optimum \citep{mladenovic1997variable}. The \textsc{Shaking} operator helps escaping local optima by searching a larger neighborhood with different basins. If this neighborhood still does not yield a better solution, an even larger neighborhood is searched until an improvement is made or the most distant neighborhood $\mathcal{N}_\kappa$ is reached. 

To search larger neighborhoods, this work proposes a destroy-and-repair operator in \textsc{Shaking} (see Algorithm~\ref{alg:shake}). During the destroy phase, $\kappa$ batches are selected and destroyed (see line 4). The selection of a batch for destruction happens through Algorithm \ref{alg:eliminate}, which is described later.

The repair phase in Algorithm~\ref{alg:shake} then reassigns the orders in the destroyed batches to the remaining batches where the additional weight opposed by the new order cannot violate the capacity constraints (see line 8). In a case where an order in a destroyed batch cannot be reassigned to any existing batches due to the capacity constraint, it will be put in a new batch $b_{|\mathcal{B}_d|+1}$ (see line 9).  Otherwise, it will be assigned to a random batch that it still fits into (see lines 11--12).

Every time an order is added to a batch, the \textsc{GreedyTour} algorithm is applied in order to assign shelves to items and update the tour in a greedy manner (see line 14).  

If the number of batches at any depot is higher than the minimum possible number of batches $|\mathcal{B}_d|^*$ (see Equation (\ref{eq:defBr})), then we will try to reduce the total number of batches through \textsc{ReduceBatches} in Algorithm \ref{alg:min_batches}, which was originally proposed by \cite{Alvim99localsearch} (see lines 16--17 of Algorithm~\ref{alg:shake}). 

\begin{algorithm}
	\SetAlgoLined
	\DontPrintSemicolon
	\caption{Eliminate a random batch}
	\label{alg:eliminate}
	\textbf{Function} \textsc{EliminateBatch}($\mathcal{B}$)\;
	\uIf{$\exists \, d \in \mathcal{V}^{\mathrm{D}}: |\mathcal{B}_d| - |\mathcal{B}_j| > 0 \, \, \, \forall \, j \in \mathcal{V}^{\mathrm{D}} \setminus d$}{
		$d \gets \underset{d\in \mathcal{V}^{\mathrm{D}}}{\mathrm{argmax}}(|\mathcal{B}_d|)$ \;
		$\mathcal{B}^{\mathrm{cand}} \gets \mathcal{B}_d $ \;
	}
	\Else{
		$\mathcal{B}^{\mathrm{cand}} \gets \mathcal{B}_d $ \;
	}
	$p_b \gets \frac{1 - (w_b / c)}{\sum_{\ell \in \mathcal{B}^{\mathrm{cand}}}1 - (w_{\ell} / c)} \quad \forall b \in \mathcal{B}^{\mathrm{cand}}$\;
	$d \gets \text{pick random } b \in \mathcal{B}^{\mathrm{cand}} \text{ according to } p_b$\;
	\textbf{return} $d$ 
\end{algorithm}
During the destroy phase, Algorithm \ref{alg:eliminate} determines whether there exists a depot $d$ with more batches than any other depot (see line 2). If this is the case, a batch belonging to this depot has to be destroyed in order not to violate the depot balancing constraints. Thus, the candidate set of batches $\mathcal{B}^{\mathrm{cand}}$ is reduced to the batches belonging to depot $d$ (line 3--4). If the depots are already balanced, all batches might be destroyed (line 6).    
The probability of a batch $b \in \mathcal{B}^{\mathrm{cand}}$ being selected is dependent on its workload, so that batches that are already close to the capacity limit have a lower chance of being selected (see line 8). After selecting a batch, the \textsc{EliminateBatch} function returns the respective batch to the \textsc{Shaking} operator, which then excludes it from the set of batches (see line 5 of Algorithm~\ref{alg:shake}).

In Algorithm~\ref{alg:min_batches}, just like in the \textsc{Shaking} operator, the \textsc{EliminateBatch} function chooses a batch $d$ to make sure that the balancing constraint is not violated (see line 3). The orders in batch $d$ are then split up among the other batches in $\mathcal{B}'$ so that the total weight of the other batches remains as small as possible; however, the capacity constraint may be violated (see lines 4--5). If no batch exceeds the maximum allowed weight, the solution $\mathcal{B}'$ is accepted (see lines 7--8); otherwise, a \textsc{Repair} operator is executed (see line 10). If the final solution still has batches that exceed the maximum weight, the neighborhood solution is discarded, otherwise it is accepted (see lines 11--12 of Algorithm \ref{alg:min_batches}).

\begin{algorithm}[h]
	\SetAlgoLined
	\DontPrintSemicolon
	\textbf{Function} \textsc{ReduceBatches}($\mathcal{B}$)\;
	\SetArgSty{textnormal}
	\While{improvement}{
		improvement $\gets$ False; $d \gets \textsc{EliminateBatch}(\mathcal{B})$; $\mathcal{B}^\prime \gets \mathcal{B} \setminus d$\;
		\For{$o \in \mathcal{O}_{d}$}{
			$b \gets \underset{b \in \mathcal{B}^\prime}{\mathrm{argmin}}(w_o + w_b)$; $\mathcal{O}_b \gets \mathcal{O}_b \cup o$\;
		}
		\uIf{$w_b \leq c \ \forall b \in \mathcal{B}^\prime$}{
			$\mathcal{B} \gets \mathcal{B}^\prime$; improvement $\gets$ True \;
		}
		\Else{
			$\bar{\mathcal{B}} \gets \textsc{Repair}(\mathcal{B}^\prime)$\;
			\uIf{$w_b \leq c \ \forall b \in \bar{\mathcal{B}}$}{
				$\mathcal{B} \gets \bar{\mathcal{B}}$; improvement $\gets$ True 
			}
		}
	}
	\textbf{return} $\mathcal{B}$ 
	\caption{Minimize number of batches}
	\label{alg:min_batches}
\end{algorithm}
\begin{algorithm}[H]
	\SetAlgoLined
	\DontPrintSemicolon
	\textbf{Function} \textsc{Repair}($\mathcal{B}^\prime$)\;
	\SetArgSty{textnormal}
	\While{\textbf{not} $w_b \leq c \ \forall b \in \mathcal{B}^\prime$ \textbf{and} improvement}{
		improvement $\gets$ False\;
		$b^{\mathrm{inf}} \gets$ draw next batch from $\{b \in \mathcal{B}\, | \, w_b > c\}$\;
		\For{$b^{\mathrm{other}} \in \mathcal{B}^\prime \setminus b^{\mathrm{inf}}$}{
			$\bar{b}^{\mathrm{inf}}, \, \bar{b}^{\mathrm{other}} \gets$ apply differencing method to $b^{\mathrm{other}}$ and $b^{\mathrm{inf}}$\;
			\If{$|\bar{b}^{\mathrm{inf}} - \bar{b}^{\mathrm{other}}| < | b^{\mathrm{inf}} - b^{\mathrm{other}} |$}{
				$\mathcal{B}^\prime \gets (\mathcal{B}^\prime \cup \{\bar{b}^{\mathrm{inf}},\bar{b}^{\mathrm{other}}\}) \setminus \{b^{\mathrm{inf}},b^{\mathrm{other}}\} $\;
				improvement $\gets$ True
			} 
		}
	}
	\textbf{return} $\mathcal{B}^\prime$ 
	\caption{Repair an infeasible neighborhood solution}
	\label{alg:repair}
\end{algorithm}




Within the \textsc{Repair} operator in Algorithm \ref{alg:repair}, the infeasible batch $b^{\mathrm{inf}}$ is fixed while iterating over all other batches to find a batch $b^{\mathrm{other}}$ to swap. The weight of $b^{\mathrm{other}}$ is as close as possible to the weight of the infeasible batch $b^{\mathrm{inf}}$ (see line 6). The reassignment of orders to $\mathcal{B}'$ is performed using the differencing method proposed by \cite{karmarkar1982differencing}. This algorithm saves the weights of all the items in a sorted list and in each iteration extracts the two largest numbers from it, computes the difference between them and places only the difference back in the list. This represents a decision to put each of the respective orders into different batches. This procedure is repeated until there is only one number left, which is the value of the final subset difference \citep{karmarkar1982differencing}.

The solution obtained by redistributing the orders is accepted if the difference in the weight of the two corresponding batches has decreased (see line 7). And the repair procedure is repeated for each batch that is infeasible or has become infeasible by applying the differencing method until all batches meet the capacity constraints or a batch cannot be balanced with any other batch. 

\subsection{Local search} \label{subsec:ls}

\begin{algorithm}[h]
	\SetAlgoLined
	\DontPrintSemicolon
	\SetArgSty{textnormal}
	\KwResult{Locally improved solution $x^*$}
	\textbf{Function} \textsc{ALS}($x, \, T, \, \Omega$)\;
	$\gamma \gets 0$; $\gamma_{\mathrm{unchanged}} \gets 0$\;
	$\rho_i \gets 1, \quad i = 1,...,|\Omega|$\tcp*{Initial weight of each operator $i$}
	\Repeat{$T \leq T^{\mathrm{min}}$ \textbf{or} $\gamma \geq \gamma^{\mathrm{max}}$ \textbf{or} $\gamma_{\mathrm{unchanged}} \geq \gamma^{\mathrm{max}}_{\mathrm{unchanged}}$}{
	    select operator $i$ using $\rho_i$\;
 	    $x^\prime \gets \Omega_i(x)$ \;
		\uIf{$U(0,1) < e^{\frac{-|x-x^\prime|}{T}}$ \textbf{or} $f(x^\prime) < f(x)$}{
		    $\gamma_{\mathrm{unchanged}} \gets 0$\;
			$x \gets x^\prime$ \tcp*{accept solution}
			\If{$f(x^\prime) < f(x^*)$}{
			    $x^* \gets x^\prime$\;
			}
		}
		\Else{
			$\gamma_{\mathrm{unchanged}} \gets \gamma_{\mathrm{unchanged}} + 1$ \tcp*{discard changes}
		}
		$T \gets \alpha T$;  $\gamma \gets \gamma + 1$\;
		$\rho_i \gets \lambda \rho_i + (1-\lambda) \Psi$ \tcp*{update the weight of $i$}
	}
	\textbf{return} $x^*$ 
	\caption{Adaptive local search}
	\label{alg:ALNS}
\end{algorithm}

After moving to a different basin of attraction using the \textsc{Shaking} algorithm, we use an adaptive local search (ALS) as described in Algorithm \ref{alg:ALNS}. We have a set of local search operators $\Omega$, and each operator $i$ has an initial weight $\rho_i = 1$ (line 3). The weights are used to select the operator using a roulette wheel principle (line 5). The algorithm calculates the probability $\phi_i$ of choosing the operator $i$ as follows: 
\begin{align}
    \phi_i = \frac{\rho_i}{\sum_{j=1}^{|\Omega|} \rho_j} \,
\end{align}
We can get a new candidate solution $x^\prime$ using the selected operator $i$ (line 6). We use the Metropolis acceptance criterion, as in simulated annealing, to accept a worse solution based on the temperature $T$ (lines 7--17). The acceptance of the worse solutions decreases if the temperature decreases. The weight of selected operator $i$ is then adjusted dynamically at the search progresses based on the quality of the new solution $x^\prime$ using the formula $\rho_i = \lambda \rho_i + (1-\lambda) \Psi$  (line 17). A score $\Psi$ is computed using the formula: 
\begin{align}
\label{eq:update}
    \Psi  =
    \mathrm{max}\begin{cases}
      \psi_1 & \text{if $f(x^\prime) < f(x^*)$}\\
      \psi_2 & \text{if $f(x^\prime) < f(x)$}\\
      \psi_3 & \text{if $f(x^\prime)$ is accepted}\\
      \psi_4 & \text{if $f(x^\prime)$ is not accepted,} \,
    \end{cases}
\end{align}
where $\psi_1, \psi_2, \psi_3,$ and $ \psi_4$ are parameters with $\psi_1 \geq \psi_2 \geq \psi_3 \geq \psi_4 \geq 0$. A high $\Psi$ value corresponds to a better performance. $\lambda \in [0,1]$ is the decay parameter that controls how sensitive the weights are to changes in the performance of the operator. Note that only the weight $\phi_i$ of the selected operator $i$ needs to be updated; the weights of other operators remain unchanged. This way of dynamically selecting an operator is adopted from adaptive variable neighborhood search (\cite{gendreau2010handbook}).

The greedy initial solution and the \textsc{Shaking} operator may come up with batches that have to visit many different shelves because each batch's items do not share the same set of storage locations. Therefore, the local search operators focus on moving orders to batches, where the number of visited shelves is minimized. In the following, we will introduce the local search operators $\Omega^{\mathrm{orders}}$ to be used in the adaptive local search. 

\begin{algorithm}[H]
	\SetAlgoLined
	\DontPrintSemicolon
	\KwResult{Updated solution $x^{\prime}$}
	\textbf{Function} \textsc{MoveOrders}($x$, $\kappa$)\;
	\SetKwFor{RepTimes}{repeat}{times}{end}
	$\mathcal{B}^{\mathrm{temp}} \gets \emptyset$\;
	determine for each order the number of required shelves to be visited: $n_o$ \;
	randomly select $\kappa$ orders according to probabilities $\phi_o = \frac{n_o}{\sum_{j \in \mathcal{O}}n_j}$   \;
	\ForEach{$i=1,...,\kappa$}{
	    $\mathcal{B}^{\mathrm{cand}} \gets \{b \in \mathcal{B} \, | \, w_b + w_{o_i} \leq c\}$ \tcp*{ensure capacity}
	    determine for each $b \in \mathcal{B}^{\mathrm{cand}}$ the number of items in $o_i$ that are stored in a shelf of $tour_b$: $\rho_b$ \;
	    randomly select $b^{\mathrm{new}}$ according to probabilities $\phi_b = \frac{\rho_b}{\sum_{j \in \mathcal{B}^{\mathrm{cand}}}\rho_j}$\;
	    $\mathcal{O}_{b^{\mathrm{old}}} \gets \mathcal{O}_{b^{\mathrm{old}}} \setminus o_i$ \;
		$\mathcal{O}_{b^{\mathrm{new}}} \gets \mathcal{O}_{b^{\mathrm{new}}} \cup o_i$ \tcp*{relocate order}
		$\mathcal{B}^{\mathrm{temp}} \gets \mathcal{B}^{\mathrm{temp}} \cup \{b^{\mathrm{old}}, b^{\mathrm{new}}\}$\;
	}
	$x^{\prime}_b \gets \textsc{OptimizeShelves}(b) \quad \forall b \in \mathcal{B}^{\mathrm{temp}}$\;
    $x^{\prime \prime} \gets \textsc{RepairNegativeStock}(x^{\prime})$\;
	\textbf{return} $x^{\prime \prime}$  
	\caption{Move $\kappa$ orders}
	\label{alg:swap_order}
\end{algorithm}

The first operator is called \textsc{MoveOrders} (see Algorithm \ref{alg:swap_order}). This operator first determines for each order $o$ the number of shelves it adds to its current batch, $n_o$, that is, the maximum number of nodes that could be eliminated by removing that order from its batch. These quantities are then normalized and used as probabilities to draw $\kappa$ orders from the set of orders $\mathcal{O}$ (lines 3--4). Given an order $o_i$, all batches that the order can still fit into are gathered in $\mathcal{B}^{\mathrm{cand}}$ (line 6). Similarly to lines 3 and 4, we determine for each of these batches the number of items of order $o_i$ that are stored in a shelf which is already approached in the tour for that batch. Again, these quantities are normalized and used as probabilities for selecting the batch that the order $o_i$ will be moved to (lines 7--8). After a new batch is selected, the order is removed from the old batch and inserted into the new one (lines 9--11). At the end of this algorithm, the tours for the changed batches will be optimized using \textsc{OptimizeShelves}, which will be explained later. Furthermore, due to the mixed-shelves storage, the tours for the batches that might share same shelves with the changed batches will also be optimized. In lines 13--15, such batches are added to $\mathcal{B}^{\mathrm{temp}}$ if they share an SKU with $o_i$. Note that lines 13--15 can be ignored for the dedicated storage.

\begin{algorithm}[H]
	\SetAlgoLined
	\DontPrintSemicolon
	\KwResult{Updated solution $x^{\prime}$}
	\textbf{Function} \textsc{ExchangeOrders}($x$, $\kappa$)\;
	$\mathcal{O}_{ij}^{\mathrm{cand}} \gets \{ $\;
	$\qquad i,j \in {\mathcal{O} \choose 2} \, | \, w_{b_{o_i}} - w_{o_i} + {w_{o_j}} \leq c \, \wedge \, w_{b_{o_j}} - w_{o_j} + {w_{o_i}} \leq c \, \wedge \, b_{o_{i}} \neq b_{o_{j}} $\;
	$\} $ \;
	determine for each order the number of shelves it adds to its batch: $\rho_o$ \;
	$\phi_{ij} \gets \frac{\rho_{o_i} + \rho_{o_j}}{\sum_{kl \in \mathcal{O}_{ij}^{\mathrm{cand}}} \rho_{o_k} + \rho_{o_l}} \quad \forall i,j \in \mathcal{O}_{ij}^{\mathrm{cand}} $\;
	$\mathcal{B}^{\mathrm{temp}} \gets \emptyset$\;
	\SetKwFor{RepTimes}{repeat}{times}{end}
	\RepTimes{$\kappa$}{
		draw $\{i,j\} \in \mathcal{O}_{ij}^{\mathrm{cand}}$ according to probabilities $\phi_{ij}  $ \;
		$\mathcal{O}_{b_{o_{i}}} \gets (\mathcal{O}_{b_{o_{i}}} \cup o_j) \setminus o_i$ \; 
		$\mathcal{O}_{b_{o_{j}}} \gets (\mathcal{O}_{b_{o_{j}}} \cup o_i) \setminus o_j$  \tcp*{exchange orders}
		$\mathcal{B}^{\mathrm{temp}} \gets \mathcal{B}^{\mathrm{temp}} \cup \{b_{o_{i}}, b_{o_{j}}\}$\;
	}
	$x^{\prime}_b \gets \textsc{OptimizeShelves}(b) \quad \forall b \in \mathcal{B}^{\mathrm{temp}}$\;
    $x^{\prime \prime} \gets \textsc{RepairNegativeStock}(x^{\prime})$\;
	\textbf{return} $x^{\prime \prime}$ 
	\caption{Exchange $\kappa$ order pairs}
	\label{alg:exchange_orders}
\end{algorithm}

The second local search operator is called \textsc{ExchangeOrders} (see Algorithm \ref{alg:exchange_orders}). At the time when the batches are highly loaded, the relocation of orders may not work any longer due to the capacity constraints. 
Therefore, the operator finds two orders $o_i$ and $o_j$ from two different batches $b_{o_{i}}$ and $b_{o_{j}}$ and exchanges them if no capacity constraints are violated. Thus the first step is to determine all combinations of orders that can be exchanged without violating any capacity constraints and come from different batches (line 2--4). Similarly to the \textsc{MoveOrders} operator, the \textsc{ExchangeOrders} operator then determines for all orders the weight $\rho_o$ specifying the number of shelves they add to their current batch (line 5). Now, for each of the possible combinations $i,j$, the weights $\rho_{i}$ and $\rho_{j}$ of the respective orders are added up and normalized in order to receive the sampling probability $\phi_{ij}$ (line 6). 
Now, in total $\kappa$ pairs of orders are selected randomly with probability $\phi_{ij}$, removed from their current batch and added to the batch of the other order in each case (lines 7--13). Similarly to \textsc{MoveOrders},  \textsc{OptimizeShelves} is used to optimize the tours for the batches in $\mathcal{B}^{\mathrm{temp}}$.

\begin{algorithm}[H]
	\SetAlgoLined
	\DontPrintSemicolon
	\KwResult{Updated solution $x^{\prime}$}
	\textbf{Function} \textsc{OptimizeShelves}($b$)\;
	$\Omega^{\mathrm{tour}} \gets \{\textsc{2-opt}, \textsc{Swap}, \textsc{Relocate}\} \,$; $T \gets \mathrm{max}(d_{ij})-\mathrm{min}(d_{ij})$\;
	$\mathcal{P}_b = \bigcup_{o \in \mathcal{O}_b} \mathcal{P}_o$; \,\, $tour_b \gets \emptyset$ \tcp*{all items of batch}	
	\While{$\mathcal{P}_b \neq \emptyset$}{
		$\mathcal{P}_{bs} \gets \{p \in \mathcal{P}_b \, | \, s \in \mathcal{V}^{\mathrm{SC}}_p \wedge n_{ps} > -9999\} \,\,\, \forall s \in \mathcal{V}^{\mathrm{SC}}$ \; 
		randomly select $s^*$ according to probabilities $\phi_s = \frac{|\mathcal{P}_{bs}|}{\sum_{j\in \mathcal{V}^{\mathrm{SC}}} |\mathcal{P}_{bj}|}$\;
		$tour_b \gets tour_b \cup s^*$\;
		\For{$p \in \mathcal{P}_{bs^*}$}{
	        $\mathcal{P}_b \gets \mathcal{P}_b \setminus p$\;
	        $n_{ps^*} \gets n_{ps^*} - 1$
		    
		}
    }
    $tour^{\prime}_b \gets \textsc{ALS}(tour_b, T, \Omega^{\mathrm{tour}})$\;
    \textbf{return} $tour^{\prime}_b$ 
\caption{Optimization of shelf assignment}
\label{alg:shelves}
\end{algorithm}

By performing the \textsc{MoveOrders} and \textsc{ExchangeOrders} operators, the new tours for the changed batches should be optimized using \textsc{OptimizeShelves} (Algorithm \ref{alg:shelves}). For a given batch, this algorithm first determines the set of items contained in it (line 3). Then, for a given batch $b$ and for each shelf, the set of items contained in batch $b$ that are stored in the respective shelf is determined (see line 5). In order to prefer solutions where potentially fewer shelves have to be visited, the cardinality of these sets $|\mathcal{P}_{bs}|$ are normalized and used as probabilities for selecting the next shelf $s^*$ to be added to the tour for $b$ (see lines 6--7). All the items that can be collected at shelf $s^*$ are then excluded from the set of items belonging to $b$, and one unit of the respective item is taken from the shelf (lines 8--11). The procedure repeats until all items of batch $b$ have been assigned to shelves.
Lastly, the resulting tour is optimized by using the ALS from algorithm \ref{alg:ALNS} and a set of local search operators $\Omega^{\mathrm{tour}}$ (line 13). Note that we use here the following local search operators. First, the \textsc{2-opt} operator exchanges two edges of a given tour to generate a neighborhood solution (see algorithm \ref{alg:2opt}). Likewise, the \textsc{Swap} operator performs the same operation, but on nodes (shelves) instead of edges connecting shelves (algorithm \ref{alg:swap}). Lastly, the \textsc{Relocate} operator moves a shelf to a new, random position in the tour (algorithm \ref{alg:relocate}).

\begin{algorithm}[H]
	\SetAlgoLined
	\DontPrintSemicolon
	\SetArgSty{textnormal}
	\textbf{Function} \textsc{2-opt}($tour := [s_1,...,s_n, s_1]$)\;
    select $i\in \{1,...,n-2\}$ and $j \in \{i+2,...,n\}$\;
	$tour^{\prime} \gets [s_1,...,s_i, s_{j-1},...,s_{i+1}, s_{j+1},...,s_{n}, s_1]$\;
	\textbf{return} $tour^\prime$ 
	\caption{2-opt}
	\label{alg:2opt}
\end{algorithm}

\begin{algorithm}[H]
	\SetAlgoLined
	\DontPrintSemicolon
	\SetArgSty{textnormal}
	\textbf{Function} \textsc{Swap}($tour := [s_1,...,s_n, s_1]$)\;
    select $ \{s_i,s_j\}  \subseteq tour$\;
	$s_i \gets s_j$, $s_j \gets s_i$ \;
	\textbf{return} $tour^\prime$ 
	\caption{Swap}
	\label{alg:swap}
\end{algorithm}

\begin{algorithm}[H]
	\SetAlgoLined
	\DontPrintSemicolon
	\SetArgSty{textnormal}
	\textbf{Function} \textsc{Relocate}($tour := [s_1,...,s_n, s_1]$)\;
    select $ \{i,j \subseteq tour\}$\;
    $tour^\prime \gets [s_1,...,s_{i-1}, s_{i+1},...,s_{j}, s_{i}, s_{j+1},...,s_{n}, s_1]$\;
	\textbf{return} $tour^\prime$ 
	\caption{Relocate}
	\label{alg:relocate}
\end{algorithm}

\begin{algorithm}[H]
	\SetAlgoLined
	\DontPrintSemicolon
	\KwResult{Feasible solution $x^{\prime}$}
	\textbf{Function} \textsc{RepairNegativeStock}($x$)\;
	\While{$\exists n_{ps} < 0 $}{
	    select infeasible item-shelf combination $(p,s)$ \;
	    $n_{ps} \gets - \infty$ \;
	    \For{each batch $b$ picking up $p$ at $s$}{
	        $x^{\mathrm{feas}}_b \gets$ \textsc{OptimizeShelves}($b$) \;
	        $\Delta_b \gets x^{\mathrm{feas}}_b - x_b$ \;
	    }
	    Select $b$ for which $\Delta_b$ is minimal; $x_b \gets x^{\mathrm{feas}}_b$; reset $n_{ps}$ \;
	}
	$x^\prime \gets \{x_b; b \in \mathcal{B}\} $ \;
    \textbf{return} $x^\prime$ 
\caption{Repair an infeasible Solution}
\label{alg:repair}
\end{algorithm}

In order to avoid the necessity of optimizing the item-shelf assignment for all batches simultaneously, we also allow for negative stocks during the assignment phase (see line 5 of algorithm \ref{alg:shelves}). Therefore, potentially infeasible solutions have to be repaired, which is done using the function \textsc{RepairNegativeStock} after applying the local search operators \textsc{MoveOrders} and \textsc{ExchangeOrders}. For a item-shelf combination $(p,s)$ with negative stock, $n_{ps}$ is set to minus infinity (line 4 in alogrithm \ref{alg:repair}), so that this combination is no longer considered in line 5 of algorithm \ref{alg:shelves}. Now, for each batch that picks up item $p$ at shelf $s$ the item-shelf assignment is applied again using \textsc{OptimizeShelves} (lines 5--8). Then, the constrained solution is accepted for the batch with the least increase in the fitness and discarded for all other batches (line 9). This procedure is repeated until there are no negative stocks anymore.
\section{Computational evaluation} \label{sec:results}
A series of numerical studies was conducted to evaluate (1) the solution qualities of three-index and two-commodity flow formulations, (2) whether the mixed-shelves policy works well for orders with different sizes of order lines (in both small and large warehouses). The computational results in this section were conducted on a PC with an Intel Core i7-7700K CPU @4.20GHz and 32GB RAM running Microsoft Windows 10 in 64-bit mode. Our models were solved by Gurobi 8.1.1 via Python version 3.7. In the following subsections, we first detail instance generation and the setting of parameters in Subsection~\ref{subsec:instance}. The results of the computational study of our models are presented in Subsection~\ref{subsec:result_model}, while the tuned parameters and the results of VNS are presented in Subsections \ref{subsec:tuning} and \ref{subsec:result_heuristics} respectively. 

\subsection{Instance generation and parameter setting} \label{subsec:instance}
\begin{table}[h]
	\centering
	\begin{tabular}{  l l l l}
		\hline
		& \multicolumn{2}{c}{Small instances}& Large instances\\
		\hline
		 Layout & $2\times12$ shelves&$6\times60$ shelves & $18\times180$ shelves\\
		 \# orders& 10, 20 & 10, 20& 100, 200\\
		 \# order lines & 1.6, 5 & 1.6, 5&1.6, 5\\
		 (average) & & & \\
	 Order sets & $a,b,c$ &$a,b,c$ &$a,b,c$\\
		 \# SKUs & 24 & 360& 3240\\
		 Storage policy & dedicated, & dedicated &dedicated,\\
		 &mixed(5)& mixed(5)& mixed(5)\\	
		 & & mixed(10) & mixed(10)\\
		 \# depots & 2 &2 & 2, 4, 8 \\	
		\# robots per depot & 1 & 1&1\\
		\hline
		 \# instance &24 & 60& 108\\
	\end{tabular}
	\captionof{table}{Parameters for instance generation.}
	\label{tab:parameters_instance_small}
\end{table}

Table~\ref{tab:parameters_instance_small} shows the parameters for small and large instances, including layout (including number of shelves), number of orders, average number of order lines, three different order sets, number of SKUs, storage policy, number of depots and number of cobots per depot. Our test instances are available at \url{https://github.com/xor-lab/rafs-datasets}.

\begin{figure}[h]
	\includegraphics[width=\textwidth]{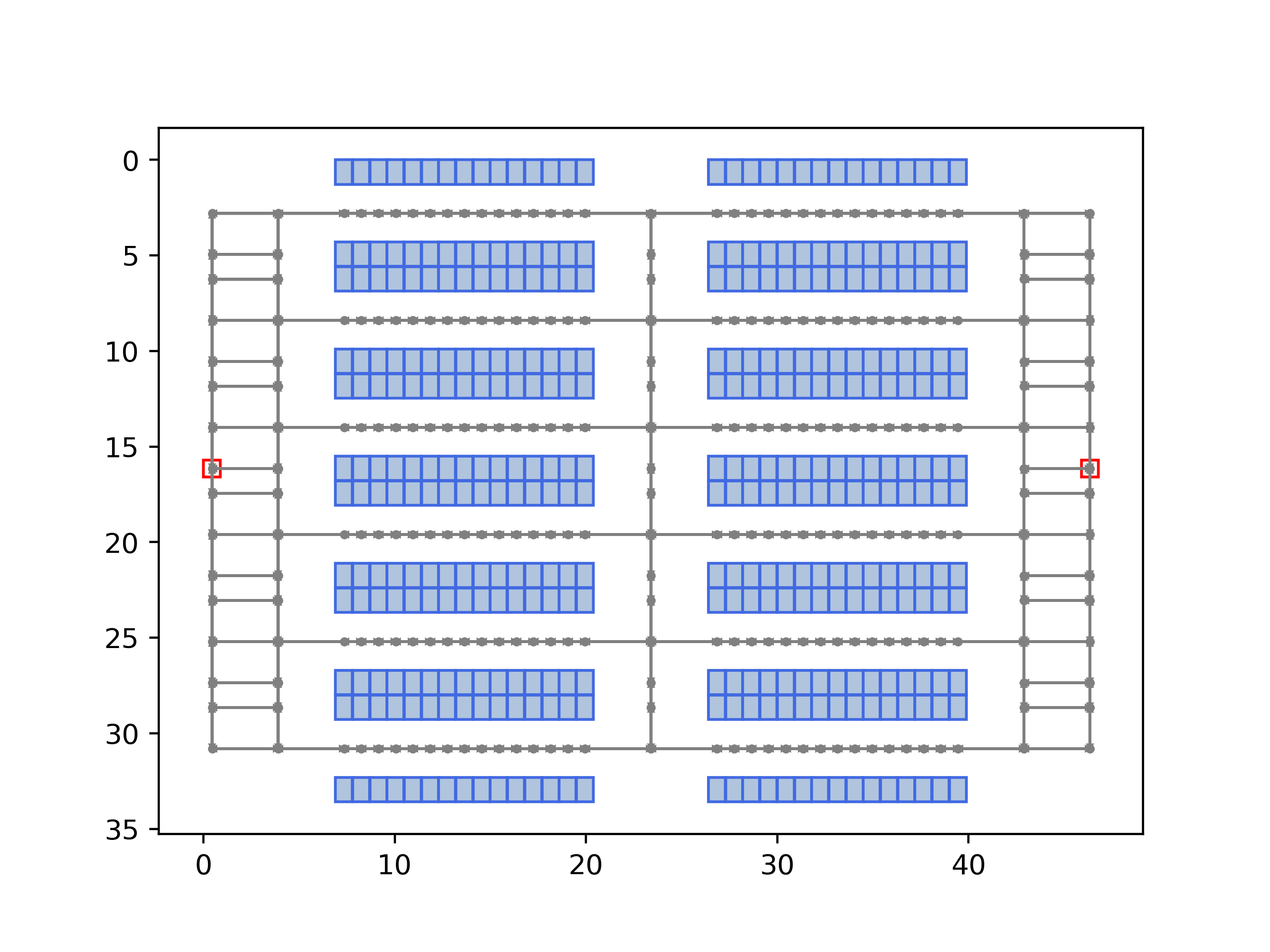}
	\caption{The layout with six pick aisles and one cross-aisle. The depots are marked in red on the left- and right-hand sides.} 
	\label{fig:layout_small}
\end{figure}

We use the same layout for traditional manual warehouses as described in \cite{van2019formulating}. For example, we have a small layout with 360 shelves (6 aisles $\times$ 60 shelves per aisle, and one cross-aisle), in Figure~\ref{fig:layout_small}. In such a layout, we define the width of each cross-aisle as 6 meters and the width of each pick aisle as 3 meters, while we define each shelf as 0.9 meters wide and 1.3 meters deep. Picking aisles are two-sided and wide enough for two-way travel of cobots (for example in \cite{lee2019robotics} for the same assumption); therefore we assume the blocking of cobots to be minimum. Another small layout we test includes 24 shelves (2 aisles $\times$ 12 shelves per aisle, without cross-aisle), while we use a large layout with 3240 shelves (18 aisles $\times$ 180 shelves per aisle, and five cross-aisles (see Figures~\ref{fig:small_layout} and \ref{fig:large_layout} in \ref{appsec:layouts})).

We generate different sets of orders (small instances: 10 and 20 orders; large instances: 100, 200 and 300 orders) randomly based on a given average number of order lines between 1 and 10. According to \cite{van2019formulating}, the number of order lines impacts the number of batches created and consequently the complexity of the problem. Therefore, we test two different average numbers of order lines per order: 1.6 (for e-commerce companies; see \cite{de2012determining}) and 5. We don't test larger average numbers of order lines than 5, since we want to keep the weight of each order smaller than the capacity of each cobot. Note that we allow more than one item to be ordered within an order line of an order. We set the capacity of each cobot as in the real-world application to be 18kg according to \cite{Wulfraat.2016}. That means the batch capacity is limited to 18kg. We call the orders with an average of 1.6 order lines \textit{small-line orders}, and the orders with an average of 5 order lines \textit{medium-line orders}.  In order to get rid of randomness by generating orders, we generate three different sets (a, b, c) for each order size. The weights of items are generated randomly between 0.5kg and 3kg (small-sized items according to practical experience).

The number of shelves sharing a SKU also seems to impact the complexity of the integrated problems (storage policy). We test two storage policies, namely dedicated and mixed-shelves. Recall that \textit{dedicated storage} means that the items of a SKU are stored in only one shelf, while \textit{mixed-shelves storage} means that the items of a SKU can be spread over several shelves. We also set the number of shelves on which the items of a SKU can be stored to 5 and 10 to see the complexity and efficiency of the mixed-shelves storage policy. We set the number of depots to be 2 for the small instances, while we set the number of depots to be 2, 4 and 8 for the large instances. For all instances, we assume that we have one cobot per depot. In our study, we do not intend to support decisions on the number of cobots and the number of depots in practice. 

Note that we get the shortest path for each two nodes $i$ and $j$ with the Dijkstra algorithm by using NetworkX (\cite{networkX}). 

\subsection{Results of mathematical models}\label{subsec:result_model}
Table~\ref{tab:gapforthreeandtwo} shows the results of the comparison between the two-commodity model (in short: two) and the three-index model (in short: three) for small instances of 24 and 360 shelves respectively. The results of applying both storage policies, namely dedicated and mixed-shelves (in short: mixed(5), mixed(10)) storage policies, for small- and medium-line orders, are also shown in the table. Note that the model for each instance was solved using eight threads. Based on the results, we can gain some knowledge of the models of integrated order batching and routing. First, the two-commodity flow formulation provides better performance by achieving smaller gaps with smaller computational time for solving our integrated problem. The improvements in the time and optimum gap are marked in green. Second, based on the results we can get the following factors that increase the complexity of the model:
\begin{itemize}
	\item the number of orders
	\item the storage policy
	\item the average number of order lines
	\item the number of shelves.
\end{itemize}
\begin{table}[h]
	\centering
	\begin{adjustbox}{width=\columnwidth,center}
	\begin{tabular}{l|l|l|r|r|r|r}
		\hline
		\multicolumn{3}{c}{}&\multicolumn{2}{|c|}{10 orders}& \multicolumn{2}{c}{20 orders} \\
		\cline{4-7}
		\multicolumn{3}{c|}{\textit{Instances with 24 shelves}}&small-line& medium-line &small-line& medium-line\\
		\hline
		\multirow{4}{*}{dedicated}& \multirow{2}{*}{three} & gap (\%)& 0&6.2&6.1&37.0\\
		& &time (s)&10.47&9212&7216&t.l.\\ 
		\cline{2-7}
		& \multirow{2}{*}{two}&gap (\%)&\cellcolor{green!25}0&\cellcolor{green!25}1.3 &\cellcolor{green!25}0&\cellcolor{green!25}22.1*\\
		& &time (s)&\cellcolor{green!25}0  &\cellcolor{green!25}5755 & \cellcolor{green!25}491.7&\cellcolor{green!25}t.l.\\
		\hline
		\multirow{4}{*}{mixed(5)}& \multirow{2}{*}{three}&gap (\%)&0&48.5&11.4& - \\
		& &time (s)&5286& t.l. & t.l. & t.l.\\
		\cline{2-7}
		& \multirow{2}{*}{two} &gap (\%)&\cellcolor{green!25}0&\cellcolor{green!25}27.8&\cellcolor{green!25}0&-\\
		& &time (s)&\cellcolor{green!25}216.3&\cellcolor{green!25}t.l. &\cellcolor{green!25}2964&t.l.\\
		\hline
		\multicolumn{3}{c|}{\textit{Instances with 360 shelves}} & \multicolumn{4}{c}{}\\
		\hline
		\multirow{4}{*}{dedicated}&\multirow{2}{*}{three}&gap (\%)&0&6.64&0&-\\
		& &time (s)&3&t.l.&2606&t.l.\\ 
		\cline{2-7}
		&\multirow{2}{*}{two}&gap (\%)&\cellcolor{green!25}0&\cellcolor{green!25}5&\cellcolor{green!25}0&-\\
		&&time(s)&\cellcolor{green!25}1.3&\cellcolor{green!25}t.l. &\cellcolor{green!25}1500&t.l.\\
		\hline
		\multirow{4}{*}{mixed(5)}& \multirow{2}{*}{three} &gap (\%)&33&-&-&\\
		& &time (s)&t.l.&t.l.&t.l.&t.l.\\
		\cline{2-7}
		& \multirow{2}{*}{two}&gap (\%)&\cellcolor{green!25}0&-&\cellcolor{green!25}28.3&-\\
		& & time (s)&\cellcolor{green!25}566.5&t.l.&\cellcolor{green!25}t.l.&t.l.\\
		\hline
		\multirow{4}{*}{mixed(10)}& \multirow{2}{*}{three} &gap (\%)&42.5&-&-&-\\
		& &time (s)&t.l.&t.l.&t.l.&t.l.\\
		\cline{2-7}
		&\multirow{2}{*}{two} &gap (\%)&\cellcolor{green!25}0&-&-&-\\
		& &time (s)&\cellcolor{green!25}2636&t.l.&t.l.&t.l.\\
	\end{tabular}
\end{adjustbox}
	\captionof{table}{Comparison of results of three-index and two-commodity models (the average gap and running time of instance sets $a,b,c$) for instances with 24 and 360 shelves. "-" means that there is no feasible integer solution found within the time limit (t.l.: 14,400 seconds)."*" means that some instances cannot find a feasible solution within the time limit.}
	\label{tab:gapforthreeandtwo}
\end{table}
\vspace{-0.3cm}
\paragraph{The number of orders} If more orders are processed in a model, then more computational time is needed. The reason is obvious, since the set of physical items $P^C$ for picking is larger for more orders, assuming all other factors are the same. 
\paragraph{The average number of order lines} The models with small-line orders are more easily solved than those with medium-line orders, since small-line orders include smaller sets of physical items $\SetOfPickingItems$ for picking (see Table~\ref{tab:gapforthreeandtwo}, columns small-line vs. medium-line), assuming the other factors are the same. Therefore the number of decision variables $x_{ijrb}$ and $z_{psrb}$ increases. Note that we should consider the number of orders and the average number of order lines together, since they influence the set of physical items for picking and the number of batches created, and consequently the complexity of the problem. Therefore, the instances of 10 medium-line orders are more difficult to solve than the instances of 20 small-line orders.
\paragraph{The storage policy} The models using the dedicated storage policy require less computational time than the models using the mixed-shelves storage policy, assuming the other factors are the same. In Table~\ref{tab:gapforthreeandtwo}, we can see that if more shelves contain items of a SKU (mixed(5) and mixed(10) vs. dedicated), the problem is more complex, since the set of possible shelves including item $p$, $V^{S}_p$, is larger for all picking items in mixed-shelves storage. Therefore the number of decision variables $x_{ijrb}$ and $z_{psrb}$ increases.
\paragraph{The number of shelves (layout)} In Table~\ref{tab:gapforthreeandtwo}, it is easy to see that the models for a larger warehouse are more difficult to solve, since the set of shelves including all items in $\SetOfPickingItems$ is larger. Therefore, the number of decision variables $x_{ijrb}$ and $z_{psrb}$ increases.

\begin{table}[h]
	\centering
	\begin{tabular}{l|r|r|r|r}
		\hline
		&\multicolumn{2}{c|}{10 orders}& \multicolumn{2}{c}{20 orders} \\
		\cline{2-5}
		&small-line & medium-line &small-line& medium-line\\
		\cline{2-5}
		\hline
		\multicolumn{5}{l}{\textit{Instances of 24 shelves}}\\
		\hline
		dedicated& 59.2 (0\%)& 220.2 (1.3\%)& 105.8 (0\%)&364.0* (22\%)\\
		\hline
		mixed(5)& 44.5 (0\%)& 139.0 (27.8\%)& 71.0 (0\%)& - \\
		\% & -24.8& -36.9& -32.9&- \\
		\hline
		\multicolumn{5}{l}{\textit{Instances of 360 shelves}}\\
		\hline
		dedicated& 232.23 (0\%)& 801.4 (9.3\%)& 352.3 (0\%)& - \\
		\hline
		mixed(5)& 110.1 (0\%)&-& 177.0 (28.3\%)& -\\
		\% & -54.1& - &-49.8 & -\\
		\hline
		mixed(10) & 89.97 (0\%)& - & -& - \\ 
		\% & -61.2 & - & - & - \\
	\end{tabular}
	\captionof{table}{The average distances for instance set $a,b,c$ of 24 and 360 shelves solved with the two-commodity formulation with the gap in \%. The lines beginning with "\%" show the percentage improvement in the results of mixed-shelves storage compared with the dedicated storage policy. "-" means that there is no feasible integer solution found within the time limit (t.l.: 14,400 seconds). "*" means that some instances cannot find a feasible solution within the time limit.}
	\label{tab:resultsmixedshelves}
\end{table}
In Table~\ref{tab:resultsmixedshelves}, the mixed-shelves policy definitely works better than the dedicated policy for small-line orders (in other words, shorter distances). In a larger warehouse (with 360 shelves), the saving is more significant. The reason is obvious: First, the distances between visiting shelves are larger in a larger warehouse. Second, it is as expected that if more shelves contain items of a SKU, more distance can be saved for picking an item of this SKU (see mixed(5) vs. mixed(10) in \textit{Instances of 360 shelves}), since the probability is higher in mixed(10) that a shelf nearby contains the required item. 

Based on the results of the two-commodity flow model in Tables~\ref{tab:gapforthreeandtwo} and \ref{tab:resultsmixedshelves}, 60\% of instances (36 out of 60) have not been solved optimally by Gurobi within the run-time limit of 4 hours. Among these, for 23 instances no feasible integer solution has been found. Layout, storage policy, the number of orders and the average number of order lines have effects on the number of non-optimal solutions. All of those instances with no feasible integer solution are those in which we apply mixed-shelves storage. This demonstrates the complexity of the problem. Therefore it is worth developing an efficient metaheuristic to gain more results. 



\subsection{Parameter tuning}\label{subsec:tuning}
It is often beneficial to tune algorithm parameters. In our VNS, parameter tuning is performed on the small instances (in total 60 instances). In our experiment, we tune two parameters in Algorithm~\ref{alg:vns}, namely $\gamma$ (parameter defining the number of iterations without improvement of $x^{\prime}$) and $\kappa_{max}$ (parameter defining the neighborhood structures). Note that we don't tune local search parameters in our experiment, such as ${\gamma}^{max}_{unchanged}=3, T^{min}=0.1$ and $\alpha=0.975$, which we apply in the local search improve routes within a batch, since they are not related to other algorithms. We test $\gamma$ with the values 100, 200, 300 and 400, while $\kappa_{max}$ is set to be 3 and 4. And we set $t_{max}$ to be 3600 seconds. Five repetitions are run to eliminate randomness. Consequently, we have in total $60\times 4\times 2 \times 5= 2400$ observations. 
Our tuning results show that the computational time increases about linearly with the increasing values of $\gamma$, but different $\gamma$-values have little effect on solution quality. So, we choose $\gamma=100$ for the rest of the experiment. From our observation, $\kappa_{max}=3$ brings similar results to the results of $\kappa_{max}=4$ but shorter computational time using the same value of $\gamma$. So we will choose $\kappa_{max}=3$ in the further experiment. 
\subsection{Results of our VNS}\label{subsec:result_heuristics}
\begin{table}[h]
	\centering
	\begin{adjustbox}{width=\columnwidth,center}
		\begin{tabular}{l|lr|lr|lr|lr}
			\hline
			&\multicolumn{4}{c|}{10 orders}& \multicolumn{4}{c}{20 orders} \\
			\cline{2-9}
			&\multicolumn{2}{c|}{small-line} & \multicolumn{2}{c|}{medium-line} &\multicolumn{2}{c|}{small-line}& \multicolumn{2}{c}{medium-line}\\
			\cline{2-9}
			\hline
			\multicolumn{9}{l}{\textit{Instances of 24 shelves}}\\
			\hline
			dedicated& \textcolor{blue}{\textbf{1}}& 59.2 (0\%)& \textcolor{blue}{\textbf{6}}&220.2 (0\%)& \textcolor{blue}{\textbf{3}}&105.8 (0\%)& \textcolor{blue}{\textbf{8}}&313.0 (302.9* (-16.8\%))\\
			\hline
			mixed(5)& \textcolor{blue}{\textbf{2}}&44.5 (0\%)& \textcolor{blue}{\textbf{7}}&132.2 (-4.9\%)& \textcolor{blue}{\textbf{5}}&71.2 (0\%)& \textcolor{blue}{\textbf{9}}&201\\
			\% & &-24.8 & &-40 & &\-32.7 & &-35.8\\
			\hline
			\multicolumn{9}{l}{\textit{Instances of 360 shelves}}\\
			\hline
			dedicated& \textcolor{blue}{\textbf{1}}&233.2 (0\%)& \textcolor{blue}{\textbf{6}}&802 (0\%)& \textcolor{blue}{\textbf{3}}&353.4 (0\%) & \textcolor{blue}{\textbf{8}}&1289 \\
			\hline
			mixed(5)& \textcolor{blue}{\textbf{2}}&110.7 (0\%)& \textcolor{blue}{\textbf{7}}&399.6 & \textcolor{blue}{\textbf{5}}&170.6 (-3.6\%)&\textcolor{blue}{\textbf{9}}& 661,2 \\
			\% && -52.5 && -50.2 && -51.7 && -48.2\\
			\hline
			mixed(10) & \textcolor{blue}{\textbf{4}}&89.97 (0\%)&\textcolor{blue}{\textbf{10}}& 304.8 & \textcolor{blue}{\textbf{11}}&153.6 & \textcolor{blue}{\textbf{12}}&545.0 \\ 
			\% & &-61.4 && -62 && -56.5 && -57.7 \\
		\end{tabular}
	\end{adjustbox}
	\captionof{table}{The average distances for instance set $a,b,c$ of 24 and 360 shelves solved with the VNS in five repetitions with parameters $\gamma=100, \kappa_{max}=3$. The lines beginning with \% shows the percentage improvement of the results of mixed-shelves storage compared with the dedicated storage policy. Each number in brackets is the percentage difference between the VNS solution and the solution solved by Gurobi in Table~\ref{tab:resultsmixedshelves}, if there is a feasible integer solution. The blue numbers in cells sort the instances based on the computational time solved with Gurobi. The number marked by * is the comparable average distance for those instances with feasible solutions solved by Gurobi.}
	\label{tab:smallresultsvns}
\end{table}
In Table~\ref{tab:smallresultsvns}, the average distances for small instances with $\kappa_{max}=3$ and $\gamma=100$ are shown. Each result represents the average value of three instance sets $a,b,c$ in five repetitions. For the instance sets where we can get optimal solutions with the two-commodity model within 4 hours, our VNS can provide optimal solutions (instance sets of 24 shelves with numbers 1, 2, 3, 5 and 6, and instance sets of 360 shelves with numbers 1 to 4, in Table~\ref{tab:smallresultsvns}). Note that the number of each instance set represents its complexity for each layout (24 and 360 shelves respectively). We want to keep the same setting with the same number in different layouts, and the instance set of 24 shelves does not include instance setting with the numbers 4, 10, 11 and 12. For the instance sets with numbers 6 to 8, some of which give us feasible solutions within 4 hours from the two-commodity model, our VNS provides either the same or better solutions. 
\begin{figure}[h]
	\begin{tikzpicture}[scale=0.85]
	\begin{axis}[xlabel={Instance setting},
	ylabel={Running time [s]},
	legend pos=north west]
	\addplot[mark=*,brown]
	coordinates
	{
		(1,0)(2, 216.3)(3,491.7)(5,2964)(6,5755)(7,14400)(8,14400)(9,14400)
	};
	\addplot[mark=o,blue]
	coordinates
	{
		(1,35)(2, 44.2)(3,37.3)(5,52.1)(6,58.4)(7,64.7)(8,85.5)(9,97.5)
	};
	\addplot[mark=*,gray]
	coordinates
	{
		(1,1.3)(2, 566.5)(3,1500)(4,2636)(5,14400)(6,14400)(7,14400)(8,14400)(9,14400)(10,14400)(11,14400)(12,14400)
	};
	\addplot[mark=o,red]
	coordinates
	{
		(1,105.2)(2, 180.0)(3,117.3)(4,241.8)(5,209)(6,118.8)(7,187.0)(8,136.7)(9,217.4)(10,247.7)(11,258.7)(12,273.1)
	};
	\legend{Model\_24,VNS\_24, Model\_360,VNS\_360}
	\end{axis}
	\end{tikzpicture}
	\caption{The comparison of running time between the two-commodity flow network model and our VNS. The number in the x-axis corresponds to the instance set with a blue number in Table~\ref{tab:smallresultsvns}.} \label{fig:timevnssolver}
\end{figure}
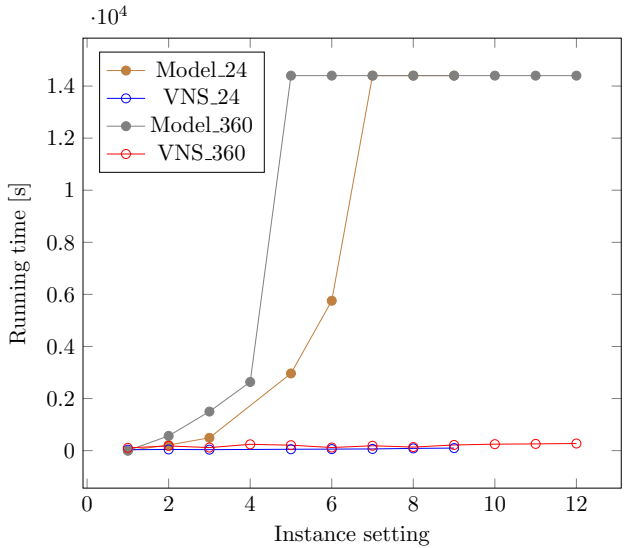

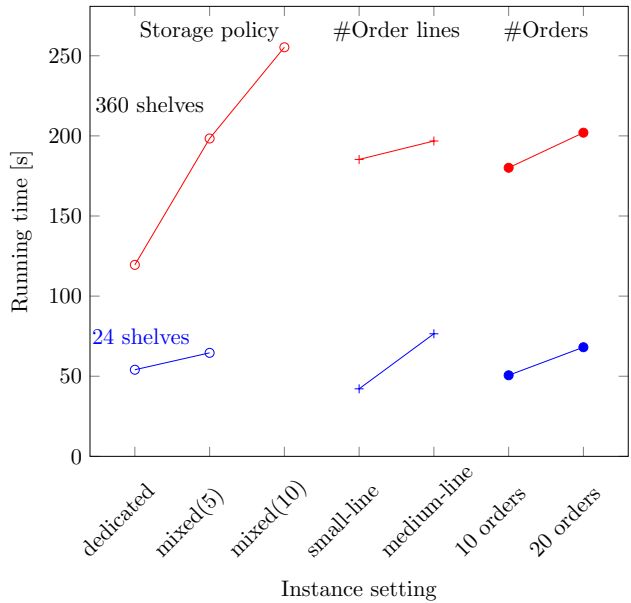
\begin{figure}[h]
	\begin{tikzpicture}[scale=0.85]
	\begin{axis}[
	legend pos=north west,
	symbolic x coords={dedicated,mixed(5),mixed(10),small-line,medium-line,10 orders,20 orders},
	xtick={dedicated,mixed(5),mixed(10),small-line,medium-line,10 orders,20 orders},
	xticklabel style={text height=0ex}, 
	ymin=0,
	xlabel=Instance setting,
	ylabel={Running time [s]},
	x tick label style={text width=2cm,rotate=45} 
	]
	\addplot[mark=o,blue]
	coordinates{
		(dedicated, 54.05)
		(mixed(5), 64.6)
	}; 
	\addplot[mark=+,blue]
	coordinates{
		(small-line, 42.15)
		(medium-line, 76.5)
	}; 
	\addplot[mark=*,blue]
	coordinates{
		(10 orders, 50.6)
		(20 orders, 68.1)
	}; 
	\draw (8,75) node {\textcolor{blue}{24 shelves}};
	\draw (100,265) node {\textcolor{black}{Storage policy}};
	\draw (350,265) node {\textcolor{black}{\#Order lines}};
	\draw (550,265) node {\textcolor{black}{\#Orders}};
	\addplot[mark=o,red]
	coordinates{
		(dedicated, 119.5)
		(mixed(5), 198.35)
		(mixed(10), 255.3)
	}; 
	\addplot[mark=+,red]
	coordinates{
		(small-line, 185.3)
		(medium-line, 196.8)
	}; 
	\addplot[mark=*,red]
	coordinates{
		(10 orders, 180.1)
		(20 orders, 202.0)
	}; 
	\draw (20,220) node {360 shelves};
	\end{axis}
	
	\end{tikzpicture}
	\caption{The average computational time of the VNS algorithm for instances which are classified by factors including the storage policy (dedicated vs. mixed(5) vs. mixed(10)), the average number of order lines (small-line vs. medium-line), the number of orders (10 orders vs. 20 orders) and the number of shelves (24 shelves vs. 360 shelves).} \label{fig:detailstimevns}
\end{figure}

The computational time of our VNS is very competitive compared with the time achieved by the solver (see Figure~\ref{fig:timevnssolver}). It decreases substantially when the problem is solved by the VNS. Note that the time is not directly comparable, since the model for each instance was solved using eight threads (to get better results for comparison), while the VNS is solved with one thread. As shown in Figure~\ref{fig:timevnssolver}, the computational time of each instance solved by the VNS is less than 5 minutes (more than 60\% of instances can be solved within 2 minutes), while most of the instances solved by the model reach the time limit (without providing optimal solutions). Figure~\ref{fig:detailstimevns} shows that the computational time applying the VNS is dependent on the layout and the storage policy. Our VNS reduces the complexity of the mixed-shelves storage (not as in the model), since the \textsc{GreedyTour} algorithm provides similar complexity for both storage policies ($\mathcal{O}(|\mathcal{P}_o|)$). The computational time increases slightly as the number of orders or order lines increases due to increased numbers of batches.  Our results show that the VNS algorithm is an efficient tool for solving the integrated order batching and routing in multi-depot mixed-shelves warehouses, at least for our small instances.


Furthermore, we can continue with our statement in Subsection~\ref{subsec:result_model} that the mixed-shelves storage policy works well for both small-line and medium-line orders. In our AGV-assisted system, the limited capacity of cobots does not allow medium-line orders with very large order lines to be accepted, unless the splitting of orders is allowed (see e.g. in \cite{xie2019efficient}). Order splitting is outside the scope of this paper. Due to this setting of limited capacity of cobots, medium-line orders in mixed-shelves storage can save up to 62\% on driving distances compared with the results in dedicated storage. The items in a medium-line order are spread among several shelves, so several shelves should be visited. Using mixed-shelves storage, the probability is higher that a shelf nearby contains the required items than dedicated storage. And because more shelves contain items of a SKU, the savings are larger for medium-line orders. The results of mixed(10) for medium-line orders are on average 11\% better than the results of mixed(5). 


\begin{table}[h]
	\centering
	\begin{adjustbox}{width=\columnwidth,center}
	\begin{tabular}{l|r|r|r|r}
		\hline
		&\multicolumn{2}{c|}{100 orders}& \multicolumn{2}{c}{200 orders} \\
		\cline{2-5}
		&small-line & medium-line &small-line& medium-line\\
		\hline
	    \multicolumn{5}{l}{\textit{2 depots}}\\
	    \hline
		dedicated& 4752.1& 18,906.0 & 9165.1&38,751.8\\
		\hline
		mixed(5)& 2895.1 & 12,657.0& 5494.5& 24,470.1\\
		\% & -39.0& -40.0 & -40.0& -42.0\\
		\hline
		mixed(10)& 2265.4 & 9715 & 4321.1& 19,664\\
		\% &-52.3 & -53.9& -52.9& -53.5\\
		\hline
	    \multicolumn{5}{l}{\textit{4 depots}}\\
	    \hline
		dedicated& 4656.8& 18,936.1 & 9397.5& 38,854.7\\
		\hline
		mixed(5)& 2843.8 &12,589.5 & 5162.6 & 23,677.4\\
		\% & -38.9 & -33.5 & -45.1& -39.1\\
		\hline
		mixed(10)& 2127.8 & 9418.6 & 3949.6 & 18,529.8\\
		\% & -54.3&-50.3 & -58.0& -52.3\\
		\hline
	    \multicolumn{5}{l}{\textit{8 depots}}\\
	    \hline
		dedicated& 4956.9& 19,623.7 & 10,053.6& 39,158.4\\
		\hline
	    mixed(5)& 2726.3 & 12,193.9 & 4978.7 & 24,283.7\\
		\% & -45.0 & -37.9& -50.5 & -38.0\\
		\hline
		mixed(10)& 1967.3 & 9107.6 & 3707.8 & 17,788.9\\
		\% & -60.3& -53.6& -63.1& -54.6\\
	\end{tabular}
	\end{adjustbox}
	\captionof{table}{The average distances for instance set $a,b,c$ of 3240 shelves (with 2, 4 and 8 depots) solved with the VNS algorithm. The lines beginning with "\%" show the percentage improvement in the results of mixed-shelves storage compared with the dedicated storage policy.}
	\label{tab:largeresultvns}
\end{table}

Saving of 54\% can also be achieved for the orders in large instances (see Table~\ref{tab:largeresultvns}, 2 depots). Based on the results of small instances, the computational time might be strongly increased for the large instances, since they have a larger layout. Therefore, we limit the computational time for solving each instance to 1 hour. Additionally, we test the large instances in the different layouts including 4 and 8 depots. The results can be different if we have other placement of the depots and shelves. Still, we can see the trend that slightly more driving distances (about 4\%) are required in the dedicated storage with the increasing size of depots. However, fewer driving distances for the instances with the small-line orders are required in the mixed-shelves storage (especially mixed(10)) with the increasing size of depots. If more shelves include items of the same SKUs, this increases the possibility that the small-line orders can be picked efficiently (regardless of in which depot). The batches are balanced among depots in large instances (same as in small instances). The mean absolute deviation of the number of assigned batches to depots is, for example, limited to less than 0.6 in the instances of 100 orders (eight depots).


\section{Conclusions} \label{sec:conclusions}

In this work, we formulate and solve integrated order batching and routing in AGV-assisted warehouses. In such warehouses, each AGV (in our case, a cobot) is sent from its depot to collect the items of orders in different locations in the storage area and goes back to its depot. Each order is collected in a bin, and that is efficient for packing. The tour for a cobot is limited by its load capacity. We consider different depots (a.k.a. loading stations for cobots) while the items for each SKU are stored in different shelves in the storage area (a.k.a. the mixed-shelves storage policy). In such an environment, we formulate the new integrated order batching and routing with three-index and two-commodity network flow formulations to minimize the driving distances for cobots. We experimentally show that the two-commodity network flow formulation works better for our integrated problem by achieving smaller gaps with smaller computational time. With mixed-shelves storage and multiple lines of orders, the complexity of the problem increases significantly. Therefore, we propose a novel variable neighborhood search (VNS) algorithm to solve the integrated problem more efficiently. According to our experiments, our VNS can obtain optimal results for small instances, which can be solved optimally in the two-commodity model within 4 hours. For those small instances where we can obtain feasible solutions from the Gurobi solver within 4 hours, our VNS algorithm provides a better solution with an improvement of up to 16.8\% within 5 minutes. Also, we apply our VNS to solve large instances within an acceptable computational time. Based on our VNS solutions, we can see that mixed-shelves storage works well for both the single-line and multi-line small-sized orders we tested. Up to 62\% can be saved on driving distances by applying the mixed-shelves storage policy compared with our reference storage policy. It is worth mentioning that the saving we achieve is under the assumption that no refill is required. In the refill process, the mixed-shelves storage may induce a larger effort when storing the items into the shelves, e.g. more shelves need to be visited to refill items of an SKU, and this may also require more shelves in the warehouse. Furthermore, we also test our large instances with different depots and we can conclude that deploying more depots does not always bring more savings of driving distances. It can bring more savings especially in cases where more shelves include items of the same SKUs when processing small-line orders. In addition, we ensure balanced batches among stations in the models and our VNS. 

Since an AGV-assisted warehouse system is a relatively new type of warehousing system, the concepts specific to such a system have not received much scholarly attention. Each picker works in a particular area in the storage area (called a \textit{zone}). There are the following limitations in this work, which can be taken as ideas for further research. First, the models proposed in this work do not consider the interaction between pickers and robots. Based on modelling the interaction in the future, the waiting time of robots for human pickers can be considered so we can extend our models to take the time as our objective and determine the batch sequencing for each cobot (in the literature, this problem is called \textit{integrated batching, routing and picker scheduling}; see e.g. \cite{van2019formulating}). Second, the waiting time of robots for human pickers in zones depends on the size of the zone; therefore it is also interesting to determine the optimal number of zones in AGV-assisted warehouse systems (see e.g. \cite{de2012determining} for the manual system). Third, we don't consider the refill process in our work. But the distribution of SKUs among the zones and within each zone affect on the one side the workload of pickers in zones, and on the other side the driving time/distance of cobots (i.e. a cobot needs more waiting time in the busy zone). Therefore, it is also interesting to research how to distribute and refill items efficiently. Mixed-shelves storage can be used as a good instrument to balance the workload of pickers, since the items of a SKU can be stored in different zones, while it reduces the waiting time for several cobots for a single picker. Fourth, the placement of depots is also an important topic for the further research.

As mentioned before, our implemented VNS is able to solve the classic integrated batching and routing problem with a single depot in the dedicated storage. As pointed out in \cite{li2017joint}, there is a lack of comparison of existing heuristics methods for solving this classic problem. Therefore, future work should be focused on find out which method in the literature performs the best.

\section*{Acknowledgements}
The authors would like to thank three anonymous referees for their insightful comments and suggestions.

\bibliography{bib.tex}

\begin{thebibliography}{41}
\expandafter\ifx\csname natexlab\endcsname\relax\def\natexlab#1{#1}\fi
\providecommand{\url}[1]{\texttt{#1}}
\providecommand{\href}[2]{#2}
\providecommand{\path}[1]{#1}
\providecommand{\DOIprefix}{doi:}
\providecommand{\ArXivprefix}{arXiv:}
\providecommand{\URLprefix}{URL: }
\providecommand{\Pubmedprefix}{pmid:}
\providecommand{\doi}[1]{\href{http://dx.doi.org/#1}{\path{#1}}}
\providecommand{\Pubmed}[1]{\href{pmid:#1}{\path{#1}}}
\providecommand{\bibinfo}[2]{#2}
\ifx\xfnm\relax \def\xfnm[#1]{\unskip,\space#1}\fi
\bibitem[{Alvim et~al.(1999)Alvim, Glover, Ribeiro \&
  Aloise}]{Alvim99localsearch}
\bibinfo{author}{Alvim, A.}, \bibinfo{author}{Glover, F.~S.},
  \bibinfo{author}{Ribeiro, C.~C.}, \& \bibinfo{author}{Aloise, D.~J.}
  (\bibinfo{year}{1999}).
\newblock \bibinfo{title}{Local search for the bin packing problem}.
\newblock In {\it \bibinfo{booktitle}{Extended Abstracts of the 3rd
  Metaheuristics International Conference}\/} (pp. \bibinfo{pages}{7--12}).
\newblock \bibinfo{publisher}{CiteSeer}.
\bibitem[{Archetti et~al.(2003)Archetti, Bertazzi \&
  Speranza}]{archetti2003reoptimizing}
\bibinfo{author}{Archetti, C.}, \bibinfo{author}{Bertazzi, L.}, \&
  \bibinfo{author}{Speranza, M.~G.} (\bibinfo{year}{2003}).
\newblock \bibinfo{title}{Reoptimizing the traveling salesman problem}.
\newblock {\it \bibinfo{journal}{Networks: An International Journal}\/},  {\it
  \bibinfo{volume}{42}\/}, \bibinfo{pages}{154--159}.
\bibitem[{Azadeh et~al.(2019)Azadeh, De~Koster \& Roy}]{azadeh2019robotized}
\bibinfo{author}{Azadeh, K.}, \bibinfo{author}{De~Koster, R.}, \&
  \bibinfo{author}{Roy, D.} (\bibinfo{year}{2019}).
\newblock \bibinfo{title}{Robotized and automated warehouse systems: Review and
  recent developments}.
\newblock {\it \bibinfo{journal}{Transportation Science}\/},  {\it
  \bibinfo{volume}{53}\/}, \bibinfo{pages}{917--945}.
  \DOIprefix\doi{10.1287/trsc.2018.0873}.
\bibitem[{Baldacci et~al.(2004)Baldacci, Hadjiconstantinou \&
  Mingozzi}]{baldacci2004exact}
\bibinfo{author}{Baldacci, R.}, \bibinfo{author}{Hadjiconstantinou, E.}, \&
  \bibinfo{author}{Mingozzi, A.} (\bibinfo{year}{2004}).
\newblock \bibinfo{title}{An exact algorithm for the capacitated vehicle
  routing problem based on a two-commodity network flow formulation}.
\newblock {\it \bibinfo{journal}{Operations Research}\/},  {\it
  \bibinfo{volume}{52}\/}, \bibinfo{pages}{723--738}.
  \DOIprefix\doi{10.1287/opre.1040.0111}.
\bibitem[{Baldacci et~al.(2007)Baldacci, Toth \& Vigo}]{baldacci2007recent}
\bibinfo{author}{Baldacci, R.}, \bibinfo{author}{Toth, P.}, \&
  \bibinfo{author}{Vigo, D.} (\bibinfo{year}{2007}).
\newblock \bibinfo{title}{Recent advances in vehicle routing exact algorithms}.
\newblock {\it \bibinfo{journal}{4OR}\/},  {\it \bibinfo{volume}{5}\/},
  \bibinfo{pages}{269--298}. \DOIprefix\doi{10.1007/s10288-007-0063-3}.
\bibitem[{Boysen et~al.(2017)Boysen, Briskorn \& Emde}]{Boysen.2017}
\bibinfo{author}{Boysen, N.}, \bibinfo{author}{Briskorn, D.}, \&
  \bibinfo{author}{Emde, S.} (\bibinfo{year}{2017}).
\newblock \bibinfo{title}{Parts-to-picker based order processing in a
  rack-moving mobile robots environment}.
\newblock {\it \bibinfo{journal}{European Journal of Operational Research}\/},
  {\it \bibinfo{volume}{262}\/}, \bibinfo{pages}{550--562}.
  \DOIprefix\doi{10.1016/j.ejor.2017.03.053}.
\bibitem[{Boysen et~al.(2019)Boysen, de~Koster \&
  Weidinger}]{boysen2019warehousing}
\bibinfo{author}{Boysen, N.}, \bibinfo{author}{de~Koster, R.}, \&
  \bibinfo{author}{Weidinger, F.} (\bibinfo{year}{2019}).
\newblock \bibinfo{title}{Warehousing in the e-commerce era: A survey}.
\newblock {\it \bibinfo{journal}{European Journal of Operational Research}\/},
  {\it \bibinfo{volume}{277}\/}, \bibinfo{pages}{396--411}.
  \DOIprefix\doi{10.1016/j.ejor.2018.08.023}.
\bibitem[{Briant et~al.(2020)Briant, Cambazard, Cattaruzza, Catusse, Ladier \&
  Ogier}]{briant2020efficient}
\bibinfo{author}{Briant, O.}, \bibinfo{author}{Cambazard, H.},
  \bibinfo{author}{Cattaruzza, D.}, \bibinfo{author}{Catusse, N.},
  \bibinfo{author}{Ladier, A.-L.}, \& \bibinfo{author}{Ogier, M.}
  (\bibinfo{year}{2020}).
\newblock \bibinfo{title}{An efficient and general approach for the joint order
  batching and picker routing problem}.
\newblock {\it \bibinfo{journal}{European Journal of Operational Research}\/},
  {\it \bibinfo{volume}{285}\/}, \bibinfo{pages}{497--512}.
  \DOIprefix\doi{10.1016/j.ejor.2020.01.059}.
\bibitem[{Chen et~al.(2015)Chen, Cheng, Chen \& Chan}]{chen2015efficient}
\bibinfo{author}{Chen, T.-L.}, \bibinfo{author}{Cheng, C.-Y.},
  \bibinfo{author}{Chen, Y.-Y.}, \& \bibinfo{author}{Chan, L.-K.}
  (\bibinfo{year}{2015}).
\newblock \bibinfo{title}{An efficient hybrid algorithm for integrated order
  batching, sequencing and routing problem}.
\newblock {\it \bibinfo{journal}{International Journal of Production
  Economics}\/},  {\it \bibinfo{volume}{159}\/}, \bibinfo{pages}{158--167}.
  \DOIprefix\doi{10.1016/j.ijpe.2014.09.029}.
\bibitem[{Cheng et~al.(2015)Cheng, Chen, Chen \& Yoo}]{cheng2015using}
\bibinfo{author}{Cheng, C.-Y.}, \bibinfo{author}{Chen, Y.-Y.},
  \bibinfo{author}{Chen, T.-L.}, \& \bibinfo{author}{Yoo, J. J.-W.}
  (\bibinfo{year}{2015}).
\newblock \bibinfo{title}{Using a hybrid approach based on the particle swarm
  optimization and ant colony optimization to solve a joint order batching and
  picker routing problem}.
\newblock {\it \bibinfo{journal}{International Journal of Production
  Economics}\/},  {\it \bibinfo{volume}{170}\/}, \bibinfo{pages}{805--814}.
  \DOIprefix\doi{10.1016/j.ijpe.2015.03.021}.
\bibitem[{Daniels et~al.(1998)Daniels, Rummel \& Schantz}]{daniels1998model}
\bibinfo{author}{Daniels, R.~L.}, \bibinfo{author}{Rummel, J.~L.}, \&
  \bibinfo{author}{Schantz, R.} (\bibinfo{year}{1998}).
\newblock \bibinfo{title}{A model for warehouse order picking}.
\newblock {\it \bibinfo{journal}{European Journal of Operational Research}\/},
  {\it \bibinfo{volume}{105}\/}, \bibinfo{pages}{1--17}.
  \DOIprefix\doi{10.1016/S0377-2217(97)00043-X}.
\bibitem[{De~Koster et~al.(2007)De~Koster, Le-Duc \& Roodbergen}]{de2007design}
\bibinfo{author}{De~Koster, R.}, \bibinfo{author}{Le-Duc, T.}, \&
  \bibinfo{author}{Roodbergen, K.~J.} (\bibinfo{year}{2007}).
\newblock \bibinfo{title}{Design and control of warehouse order picking: A
  literature review}.
\newblock {\it \bibinfo{journal}{European {J}ournal of {O}perational
  {R}esearch}\/},  {\it \bibinfo{volume}{182}\/}, \bibinfo{pages}{481--501}.
  \DOIprefix\doi{10.1016/j.ejor.2006.07.009}.
\bibitem[{De~Koster et~al.(2012)De~Koster, Le-Duc \&
  Zaerpour}]{de2012determining}
\bibinfo{author}{De~Koster, R.~B.}, \bibinfo{author}{Le-Duc, T.}, \&
  \bibinfo{author}{Zaerpour, N.} (\bibinfo{year}{2012}).
\newblock \bibinfo{title}{Determining the number of zones in a pick-and-sort
  order picking system}.
\newblock {\it \bibinfo{journal}{International Journal of Production
  Research}\/},  {\it \bibinfo{volume}{50}\/}, \bibinfo{pages}{757--771}.
  \DOIprefix\doi{10.1080/00207543.2010.543941}.
\bibitem[{Ene \& {\"O}zt{\"u}rk(2012)}]{ene2012storage}
\bibinfo{author}{Ene, S.}, \& \bibinfo{author}{{\"O}zt{\"u}rk, N.}
  (\bibinfo{year}{2012}).
\newblock \bibinfo{title}{Storage location assignment and order picking
  optimization in the automotive industry}.
\newblock {\it \bibinfo{journal}{The International Journal of Advanced
  Manufacturing Technology}\/},  {\it \bibinfo{volume}{60}\/},
  \bibinfo{pages}{787--797}. \DOIprefix\doi{10.1007/s00170-011-3593-y}.
\bibitem[{Fabric(2020)}]{fabric.2020}
\bibinfo{author}{Fabric} (\bibinfo{year}{2020}).
\newblock \bibinfo{title}{{The impact of Covid-19 on online grocery: Why the
  coronavirus pandemic in 2020 will change the grocery industry forever in the
  United States}}.
\newblock \URLprefix
  \url{https://getfabric.com/the-impact-of-covid-19-on-online-grocery/}.
\bibitem[{Finke(1984)}]{finke1984two}
\bibinfo{author}{Finke, G.} (\bibinfo{year}{1984}).
\newblock \bibinfo{title}{A two-commodity network flow approach to the
  traveling salesman problem}.
\newblock {\it \bibinfo{journal}{Congressus Numerantium}\/},  {\it
  \bibinfo{volume}{41}\/}, \bibinfo{pages}{167--178}.
\bibitem[{Garey \& Johnson(1979)}]{garey1979computers}
\bibinfo{author}{Garey, M.~R.}, \& \bibinfo{author}{Johnson, D.~S.}
  (\bibinfo{year}{1979}).
\newblock {\it \bibinfo{title}{Computers and intractability: {A} guide to the
  theory of {NP}-completeness}\/} volume \bibinfo{volume}{174}.
\newblock \bibinfo{publisher}{Freeman}.
\bibitem[{Gendreau \& Potvin(2010)}]{gendreau2010handbook}
\bibinfo{author}{Gendreau, M.}, \& \bibinfo{author}{Potvin, J.-Y.~E.}
  (\bibinfo{year}{2010}).
\newblock {\it \bibinfo{title}{Handbook of metaheuristics}\/}
  volume~\bibinfo{volume}{3}.
\newblock \bibinfo{publisher}{Springer}.
\newblock \DOIprefix\doi{10.1007/978-3-319-91086-4}.
\bibitem[{Giosa et~al.(2002)Giosa, Tansini \& Viera}]{giosa2002new}
\bibinfo{author}{Giosa, I.}, \bibinfo{author}{Tansini, I.}, \&
  \bibinfo{author}{Viera, I.} (\bibinfo{year}{2002}).
\newblock \bibinfo{title}{New assignment algorithms for the multi-depot vehicle
  routing problem}.
\newblock {\it \bibinfo{journal}{Journal of the Operational Research
  Society}\/},  {\it \bibinfo{volume}{53}\/}, \bibinfo{pages}{977--984}.
  \DOIprefix\doi{10.1057/palgrave.jors.2601426}.
\bibitem[{Hagberg et~al.(2008)Hagberg, Schult \& Swart}]{networkX}
\bibinfo{author}{Hagberg, A.~A.}, \bibinfo{author}{Schult, D.~A.}, \&
  \bibinfo{author}{Swart, P.~J.} (\bibinfo{year}{2008}).
\newblock \bibinfo{title}{Exploring network structure, dynamics, and function
  using {NetworkX}}.
\newblock In \bibinfo{editor}{G.~Varoquaux}, \bibinfo{editor}{T.~Vaught}, \&
  \bibinfo{editor}{J.~Millman} (Eds.), {\it \bibinfo{booktitle}{Proceedings of
  the 7th Python in Science Conference}\/} (pp. \bibinfo{pages}{11--15}).
\newblock \bibinfo{address}{Pasadena, CA, USA}.
\bibitem[{Karmarkar \& Karp(1982)}]{karmarkar1982differencing}
\bibinfo{author}{Karmarkar, N.}, \& \bibinfo{author}{Karp, R.~M.}
  (\bibinfo{year}{1982}).
\newblock {\it \bibinfo{title}{The differencing method of set partitioning}\/}.
\newblock \bibinfo{publisher}{Computer Science Division (EECS), University of
  California}.
\bibitem[{Kulak et~al.(2012)Kulak, Sahin \& Taner}]{kulak2012joint}
\bibinfo{author}{Kulak, O.}, \bibinfo{author}{Sahin, Y.}, \&
  \bibinfo{author}{Taner, M.~E.} (\bibinfo{year}{2012}).
\newblock \bibinfo{title}{Joint order batching and picker routing in single and
  multiple-cross-aisle warehouses using cluster-based tabu search algorithms}.
\newblock {\it \bibinfo{journal}{Flexible Services and Manufacturing
  Journal}\/},  {\it \bibinfo{volume}{24}\/}, \bibinfo{pages}{52--80}.
  \DOIprefix\doi{10.1007/s10696-011-9101-8}.
\bibitem[{Kyt{\"o}joki et~al.(2007)Kyt{\"o}joki, Nuortio, Br{\"a}ysy \&
  Gendreau}]{kytojoki2007efficient}
\bibinfo{author}{Kyt{\"o}joki, J.}, \bibinfo{author}{Nuortio, T.},
  \bibinfo{author}{Br{\"a}ysy, O.}, \& \bibinfo{author}{Gendreau, M.}
  (\bibinfo{year}{2007}).
\newblock \bibinfo{title}{An efficient variable neighborhood search heuristic
  for very large scale vehicle routing problems}.
\newblock {\it \bibinfo{journal}{Computers \& Operations Research}\/},  {\it
  \bibinfo{volume}{34}\/}, \bibinfo{pages}{2743--2757}.
  \DOIprefix\doi{10.1016/j.cor.2005.10.010}.
\bibitem[{Lee \& Murray(2019)}]{lee2019robotics}
\bibinfo{author}{Lee, H.-Y.}, \& \bibinfo{author}{Murray, C.~C.}
  (\bibinfo{year}{2019}).
\newblock \bibinfo{title}{Robotics in order picking: {E}valuating warehouse
  layouts for pick, place, and transport vehicle routing systems}.
\newblock {\it \bibinfo{journal}{International Journal of Production
  Research}\/},  {\it \bibinfo{volume}{57}\/}, \bibinfo{pages}{5821--5841}.
  \DOIprefix\doi{10.1080/00207543.2018.1552031}.
\bibitem[{Li et~al.(2017)Li, Huang \& Dai}]{li2017joint}
\bibinfo{author}{Li, J.}, \bibinfo{author}{Huang, R.}, \& \bibinfo{author}{Dai,
  J.~B.} (\bibinfo{year}{2017}).
\newblock \bibinfo{title}{Joint optimisation of order batching and picker
  routing in the online retailer’s warehouse in {China}}.
\newblock {\it \bibinfo{journal}{International Journal of Production
  Research}\/},  {\it \bibinfo{volume}{55}\/}, \bibinfo{pages}{447--461}.
  \DOIprefix\doi{10.1080/00207543.2016.1187313}.
\bibitem[{Lin et~al.(2016)Lin, Kang, Hou \& Cheng}]{lin2016joint}
\bibinfo{author}{Lin, C.-C.}, \bibinfo{author}{Kang, J.-R.},
  \bibinfo{author}{Hou, C.-C.}, \& \bibinfo{author}{Cheng, C.-Y.}
  (\bibinfo{year}{2016}).
\newblock \bibinfo{title}{Joint order batching and picker {Manhattan} routing
  problem}.
\newblock {\it \bibinfo{journal}{Computers \& Industrial Engineering}\/},  {\it
  \bibinfo{volume}{95}\/}, \bibinfo{pages}{164--174}.
  \DOIprefix\doi{10.1016/j.cie.2016.03.009}.
\bibitem[{Matusiak et~al.(2014)Matusiak, de~Koster, Kroon \&
  Saarinen}]{matusiak2014fast}
\bibinfo{author}{Matusiak, M.}, \bibinfo{author}{de~Koster, R.},
  \bibinfo{author}{Kroon, L.}, \& \bibinfo{author}{Saarinen, J.}
  (\bibinfo{year}{2014}).
\newblock \bibinfo{title}{A fast simulated annealing method for batching
  precedence-constrained customer orders in a warehouse}.
\newblock {\it \bibinfo{journal}{European Journal of Operational Research}\/},
  {\it \bibinfo{volume}{236}\/}, \bibinfo{pages}{968--977}.
  \DOIprefix\doi{10.1016/j.ejor.2013.06.001}.
\bibitem[{Mladenovi{\'c} \& Hansen(1997)}]{mladenovic1997variable}
\bibinfo{author}{Mladenovi{\'c}, N.}, \& \bibinfo{author}{Hansen, P.}
  (\bibinfo{year}{1997}).
\newblock \bibinfo{title}{Variable neighborhood search}.
\newblock {\it \bibinfo{journal}{Computers \& Operations Research}\/},  {\it
  \bibinfo{volume}{24}\/}, \bibinfo{pages}{1097--1100}.
  \DOIprefix\doi{10.1016/S0305-0548(97)00031-2}.
\bibitem[{Montoya-Torres et~al.(2015)Montoya-Torres, Franco, Isaza, Jim{\'e}nez
  \& Herazo-Padilla}]{montoya2015literature}
\bibinfo{author}{Montoya-Torres, J.~R.}, \bibinfo{author}{Franco, J.~L.},
  \bibinfo{author}{Isaza, S.~N.}, \bibinfo{author}{Jim{\'e}nez, H.~F.}, \&
  \bibinfo{author}{Herazo-Padilla, N.} (\bibinfo{year}{2015}).
\newblock \bibinfo{title}{A literature review on the vehicle routing problem
  with multiple depots}.
\newblock {\it \bibinfo{journal}{Computers \& Industrial Engineering}\/},  {\it
  \bibinfo{volume}{79}\/}, \bibinfo{pages}{115--129}.
  \DOIprefix\doi{10.1016/j.cie.2014.10.029}.
\bibitem[{Ramos et~al.(2020)Ramos, Gomes \& P{\'o}voa}]{ramos2019multi}
\bibinfo{author}{Ramos, T. R.~P.}, \bibinfo{author}{Gomes, M.~I.}, \&
  \bibinfo{author}{P{\'o}voa, A. P.~B.} (\bibinfo{year}{2020}).
\newblock \bibinfo{title}{Multi-depot vehicle routing problem: {A} comparative
  study of alternative formulations}.
\newblock {\it \bibinfo{journal}{International Journal of Logistics Research
  and Applications}\/},  {\it \bibinfo{volume}{23}\/},
  \bibinfo{pages}{103--120}. \DOIprefix\doi{10.1080/13675567.2019.1630374}.
\bibitem[{Scholz et~al.(2017)Scholz, Schubert \& W{\"a}scher}]{scholz2017order}
\bibinfo{author}{Scholz, A.}, \bibinfo{author}{Schubert, D.}, \&
  \bibinfo{author}{W{\"a}scher, G.} (\bibinfo{year}{2017}).
\newblock \bibinfo{title}{Order picking with multiple pickers and due dates:
  {S}imultaneous solution of order batching, batch assignment and sequencing,
  and picker routing problems}.
\newblock {\it \bibinfo{journal}{European Journal of Operational Research}\/},
  {\it \bibinfo{volume}{263}\/}, \bibinfo{pages}{461--478}.
  \DOIprefix\doi{10.1016/j.ejor.2017.04.038}.
\bibitem[{Tompkins(2010)}]{Tompkins.2010}
\bibinfo{author}{Tompkins, J.~A.} (\bibinfo{year}{2010}).
\newblock {\it \bibinfo{title}{Facilities planning}\/}.
\newblock (\bibinfo{edition}{4th} ed.).
\newblock \bibinfo{publisher}{{John Wiley {\&} Sons}}.
\bibitem[{Tsai et~al.(2008)Tsai, Liou \& Huang}]{tsai2008using}
\bibinfo{author}{Tsai, C.-Y.}, \bibinfo{author}{Liou, J.~J.}, \&
  \bibinfo{author}{Huang, T.-M.} (\bibinfo{year}{2008}).
\newblock \bibinfo{title}{Using a multiple-{GA} method to solve the batch
  picking problem: {C}onsidering travel distance and order due time}.
\newblock {\it \bibinfo{journal}{International Journal of Production
  Research}\/},  {\it \bibinfo{volume}{46}\/}, \bibinfo{pages}{6533--6555}.
  \DOIprefix\doi{10.1080/00207540701441947}.
\bibitem[{Valle et~al.(2017)Valle, Beasley \& da~Cunha}]{valle2017optimally}
\bibinfo{author}{Valle, C.~A.}, \bibinfo{author}{Beasley, J.~E.}, \&
  \bibinfo{author}{da~Cunha, A.~S.} (\bibinfo{year}{2017}).
\newblock \bibinfo{title}{Optimally solving the joint order batching and picker
  routing problem}.
\newblock {\it \bibinfo{journal}{European Journal of Operational Research}\/},
  {\it \bibinfo{volume}{262}\/}, \bibinfo{pages}{817--834}.
  \DOIprefix\doi{10.1016/j.ejor.2017.03.069}.
\bibitem[{Van~Gils et~al.(2019)Van~Gils, Caris, Ramaekers \&
  Braekers}]{van2019formulating}
\bibinfo{author}{Van~Gils, T.}, \bibinfo{author}{Caris, A.},
  \bibinfo{author}{Ramaekers, K.}, \& \bibinfo{author}{Braekers, K.}
  (\bibinfo{year}{2019}).
\newblock \bibinfo{title}{Formulating and solving the integrated batching,
  routing, and picker scheduling problem in a real-life spare parts warehouse}.
\newblock {\it \bibinfo{journal}{European Journal of Operational Research}\/},
  {\it \bibinfo{volume}{277}\/}, \bibinfo{pages}{814--830}.
  \DOIprefix\doi{10.1016/j.ejor.2019.03.012}.
\bibitem[{Van~Gils et~al.(2018)Van~Gils, Ramaekers, Caris \&
  De~Koster}]{van2018designing}
\bibinfo{author}{Van~Gils, T.}, \bibinfo{author}{Ramaekers, K.},
  \bibinfo{author}{Caris, A.}, \& \bibinfo{author}{De~Koster, R.~B.}
  (\bibinfo{year}{2018}).
\newblock \bibinfo{title}{Designing efficient order picking systems by
  combining planning problems: State-of-the-art classification and review}.
\newblock {\it \bibinfo{journal}{European Journal of Operational Research}\/},
  {\it \bibinfo{volume}{267}\/}, \bibinfo{pages}{1--15}.
  \DOIprefix\doi{10.1016/j.ejor.2017.09.002}.
\bibitem[{Weidinger(2018)}]{weidinger2018picker}
\bibinfo{author}{Weidinger, F.} (\bibinfo{year}{2018}).
\newblock \bibinfo{title}{Picker routing in rectangular mixed shelves
  warehouses}.
\newblock {\it \bibinfo{journal}{Computers \& Operations Research}\/},  {\it
  \bibinfo{volume}{95}\/}, \bibinfo{pages}{139--150}.
  \DOIprefix\doi{10.1016/j.cor.2018.03.012}.
\bibitem[{Weidinger et~al.(2019)Weidinger, Boysen \&
  Schneider}]{weidinger2019picker}
\bibinfo{author}{Weidinger, F.}, \bibinfo{author}{Boysen, N.}, \&
  \bibinfo{author}{Schneider, M.} (\bibinfo{year}{2019}).
\newblock \bibinfo{title}{Picker routing in the mixed-shelves warehouses of
  e-commerce retailers}.
\newblock {\it \bibinfo{journal}{European Journal of Operational Research}\/},
  {\it \bibinfo{volume}{274}\/}, \bibinfo{pages}{501--515}.
  \DOIprefix\doi{10.1016/j.ejor.2018.10.021}.
\bibitem[{Won \& Olafsson(2005)}]{won2005joint}
\bibinfo{author}{Won, J.}, \& \bibinfo{author}{Olafsson, S.}
  (\bibinfo{year}{2005}).
\newblock \bibinfo{title}{Joint order batching and order picking in warehouse
  operations}.
\newblock {\it \bibinfo{journal}{International Journal of Production
  Research}\/},  {\it \bibinfo{volume}{43}\/}, \bibinfo{pages}{1427--1442}.
  \DOIprefix\doi{10.1080/00207540410001733896}.
\bibitem[{Wulfraat(2016)}]{Wulfraat.2016}
\bibinfo{author}{Wulfraat, M.} (\bibinfo{year}{2016}).
\newblock \bibinfo{title}{{L}ocus {R}obotics: {A}n independent consultant
  review of autonomous robots}.
\newblock \URLprefix
  \url{https://mwpvl.com/html/locus_robotics_-_independent_consultant_review.html}.
\bibitem[{Xie et~al.(2020)Xie, Thieme, Krenzler \& Li}]{xie2019efficient}
\bibinfo{author}{Xie, L.}, \bibinfo{author}{Thieme, N.},
  \bibinfo{author}{Krenzler, R.}, \& \bibinfo{author}{Li, H.}
  (\bibinfo{year}{2020}).
\newblock \bibinfo{title}{Introducing split orders and optimizing operational
  policies in robotic mobile fulfillment systems}.
\newblock {\it \bibinfo{journal}{European Journal of Operational Research}\/},
  . \DOIprefix\doi{https://doi.org/10.1016/j.ejor.2020.05.032}.

\end{thebibliography}

\appendix
\section{Components of an AVG-assisted picking system} \label{subsec:rmfs}
Firstly, we define some terms related to orders in Table \ref{tab:terms}.

\begin{table}[H]
	\centering
	\begin{tabular}{ l l}
		\hline
		\textbf{Term} & \textbf{Description}\\
		\hline
		\textit{SKU} & stock keeping unit \\
		\hline
		\textit{order line} & one SKU with the ordered quantity \\
		\hline
		\textit{item} & a physical unit of one SKU\\
		\hline
		\textit{pick order}& a set of order lines from a customer's order\\
		\hline
		\textit{small-line order}& small number of SKUs (on average 1.6) per order\\
		\hline
		\textit{multi-line order}& multiple SKUs per order\\
		\hline
		\textit{backlog} & all unfulfilled orders\\
		\hline
	\end{tabular}
	\captionof{table}{Terms related to orders.}
	\label{tab:terms}
\end{table}

The central components of an AGV-assisted picking system are listed in Table~\ref{tab:components}.

\begin{table}[H]
	\centering
	\begin{tabular}{ l l}
		\hline
		\textbf{Component} & \textbf{Description}\\
		\hline
		\textit{shelves} & storing items\\
		\hline
		\textit{bin} & temporary storage items for an order during picking\\
		\hline
		\textit{storage area} & the inventory area where the shelves are stored\\
		\hline	
		human workers & \\	
		\textit{picker} & working in a given part of the storage area to pick \\
		&items from shelves \\
		\textit{packer} & packing order items in packing stations\\
		\hline
		\textit{packing stations} & where packers pack the pick order items\\
		\hline
		\textit{queue} & a list of cobots waiting in front of a packing station \\
		\hline
		\textit{picking locations} & somewhere in front of the shelves in the storage area\\
		& where pickers pick the order items \\
		\hline
		\textit{cobots} & robots that transport items between the storage area and\\
		& workstations. Each of them has a limited weight capacity \\
		& and belongs to a certain packing station.\\
		\hline
	\end{tabular}
	\captionof{table}{The central components of an AGV-assisted picking system.}
	\label{tab:components}
\end{table}

\section{Additional layouts used in our experiments}\label{appsec:layouts}
In addition to the layout with six pick aisles and one cross-aisle, we also use other layouts as shown in the following figures.
\begin{figure}[h]
    \centering
    \includegraphics[width=0.7\textwidth]{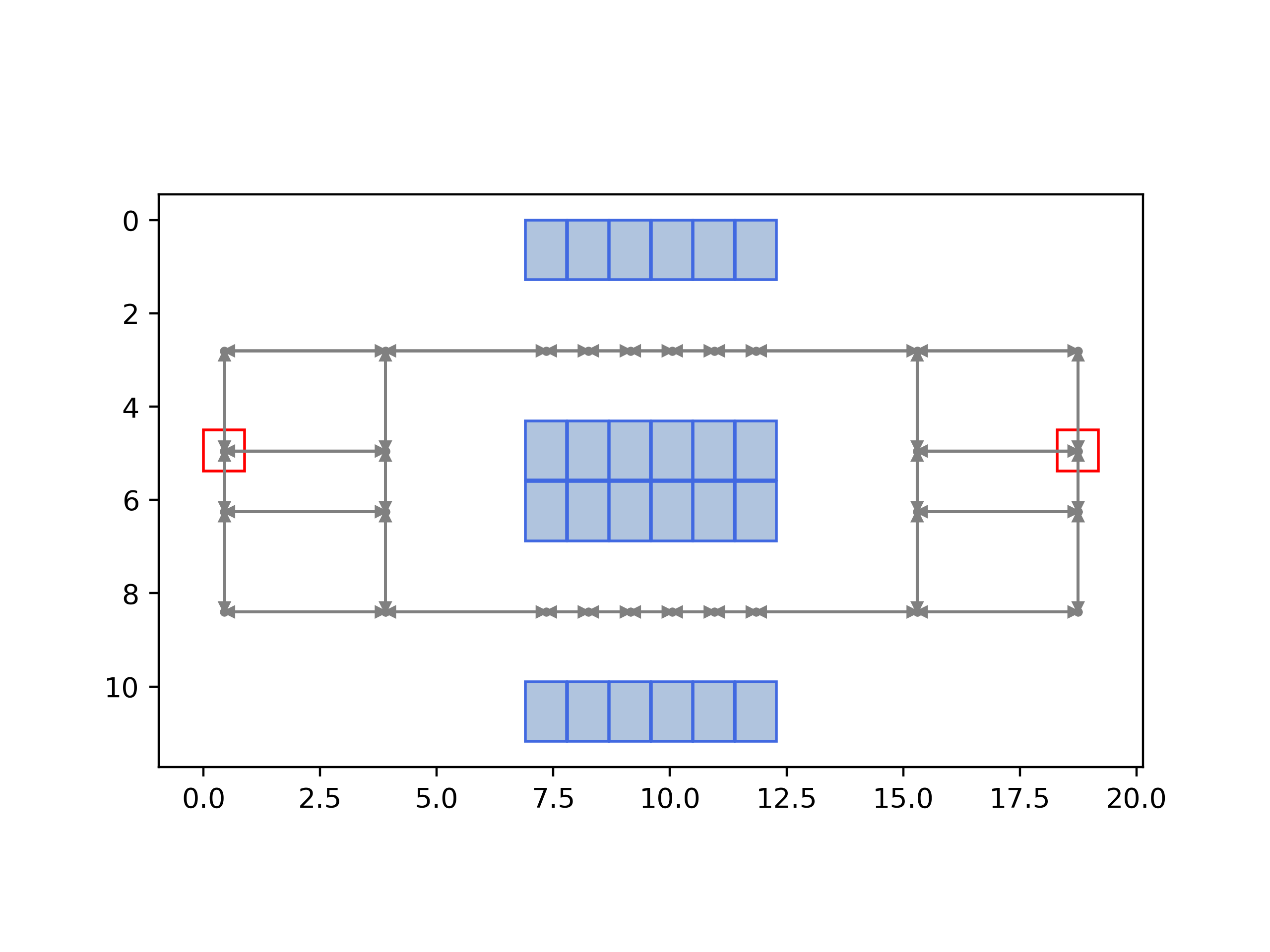}
    \caption{The small layout with two pick aisles and without cross-aisle.}
    \label{fig:small_layout}
\end{figure}
\vspace{-0.5cm}
\begin{figure}[H]
    \centering
    \includegraphics[width=\textwidth]{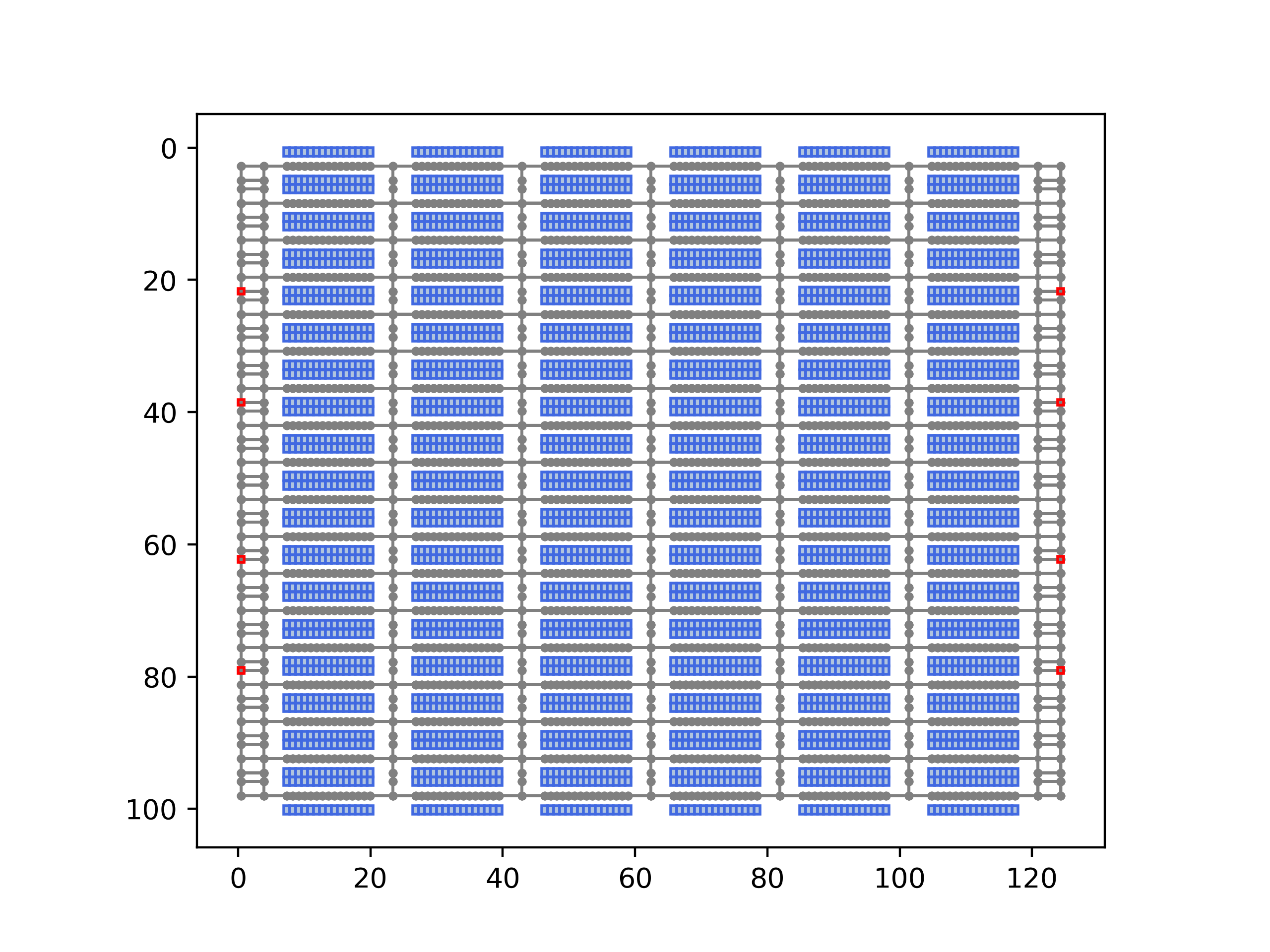}
    \caption{The large layout with 18 pick aisles, 5 cross-aisles and 8 depots (4 on the left and 4 on the right).}
    \label{fig:large_layout}
\end{figure}
\end{document}